\documentclass[pdflatex,sn-mathphys-num]{sn-jnl}% Math and Physical 
\usepackage[utf8]{inputenc}
\usepackage{amsmath}
\usepackage{amssymb}
\usepackage{amsfonts}
\usepackage{xcolor}
\usepackage{bbold}
\usepackage{bbm}
\usepackage{bm}
\usepackage{hyperref}
\usepackage{ulem}

\usepackage{algorithm}
\usepackage{algpseudocode}

\algnewcommand{\IIf}[1]{\State\algorithmicif\ #1\ \algorithmicthen}
\algnewcommand{\EndIIf}{\unskip\ \algorithmicend\ \algorithmicif}

\def\eps{\epsilon}

\newtheorem{assumption}{Assumption}
\newtheorem{defin}[assumption]{Definition}
\newtheorem{lemma}[assumption]{Lemma}
\newtheorem{corollary}[assumption]{Corollary}
\newtheorem{proposition}[assumption]{Proposition}
\newtheorem{remark}[assumption]{Remark}
\newtheorem{theorem}[assumption]{Theorem}

\title[A linesearch-based derivative-free method for noisy black-box problems]{A linesearch-based derivative-free method for noisy black-box problems}

\author[1]{\fnm{Alberto} \sur{De Santis}}\email{desantis@diag.uniroma1.it}
\equalcont{These authors contributed equally to this work.}

\author*[1,2]{\fnm{Giampaolo} \sur{Liuzzi}}\email{liuzzi@diag.uniroma1.it}
\equalcont{These authors contributed equally to this work.}

\author[1]{\fnm{Stefano} \sur{Lucidi}}\email{lucidi@diag.uniroma1.it}
\equalcont{These authors contributed equally to this work.}

\affil*[1]{\orgdiv{Department of Computer Control and Management Engineering ``A. Ruberti''}, \orgname{``Sapienza'' University of Rome}, \orgaddress{\city{Rome}, \country{Italy}}}

\affil*[2]{\orgdiv{Institute for Systems Analysis and Computer Science ``A. Ruberti''}, \orgname{National Research Council}, \orgaddress{\city{Rome}, \country{Italy}}}

\renewcommand{\Re}{\mathbb{R}}

\begin{document}

\abstract{In this work we consider unconstrained optimization problems. The objective function is known through a zeroth order stochastic oracle that gives an estimate of the true objective function. To solve these problems, we propose a derivative-free algorithm based on extrapolation techniques. Under reasonable assumptions we are able to prove convergence properties for the proposed algorithms. Furthermore, we also give a worst-case complexity result stating that the total number of iterations where the expected value of the norm of the objective function gradient is above a prefixed $\eps>0$ is ${\cal O}(n^2\eps^{-2}/\beta^2)$ in the worst case.} 


\maketitle

% {\small
% {\bf Abstract} ...
% \par\medskip

% {\bf Keywords}: ...
% \par\medskip

%{\bf MSC subject classification}: ...

%\par\medskip
%}

\nocite{dzahini2023constrained}
\nocite{rinaldi2024stochastic}
\nocite{chen2018stochastic}
\nocite{blanchet2019convergence}
\nocite{audet2021stochastic}

\section{Introduction}

In this paper we consider the following problem
\begin{equation}\label{theoretical_prob}
   \min_{x\in\Re^n} f(x)\qquad \text{where}\ f(x) = \mathbb{E}_\theta [f(x,\theta)]
\end{equation}
where $f:\Re^n\to\Re$. We assume (as it is common) that  first-order derivatives are not available or unreliable to obtain. The framework of problem \eqref{theoretical_prob} is very general in that it encompasses many relevant applications and real world problems, see e.g. \cite{bottou2018optimization,amaran2016simulation}. Derivative-free methods considering problem \eqref{theoretical_prob} have been recently proposed in the literature. They mainly belong to the following general classes of methods: model-based methods and direct search methods. 

%\cite{menickelly2023stochastic} % mi pare che usa anche un oracolo dei gradienti con rumore.

{\em Model-based methods}. In \cite{larson2016stochastic} a trust region method based on probabilistic models is proposed to minimize a stochastic function. Convergence to zero of the norm of the gradient is shown in probability. 
In \cite{shashaani2018astro} the ASTRO-DF algorithm has been proposed which tackles the stochasticity of the objective function through automatic resampling of function values. In \cite{chen2018stochastic} a trust region method is proposed which is based on the use of probabilistic models. Complexity analysis of the algorithm proposed in \cite{chen2018stochastic} is carried out in \cite{blanchet2019convergence}, via the theory of supermartingales. 
In \cite{rinaldi2024stochastic}, a trust-region algorithm is proposed along with a condition on probabilistic estimates that allows to save function samplings with respect to other available methods. 

{\em Direct search methods.} In \cite{alarie2021optimization} a Mesh Adaptive Direct Search (MADS) algorithm with dynamic precision is proposed. The algorithm and convergence analysis are based on the assumption that it is possible to control the precision of the oracle computing function value estimates.
In \cite{audet2021stochastic} the StoMADS algorithm, a MADS algorithm for stochastic problems, is proposed. StoMADS has been modified in \cite{dzahini2023constrained} to also tackle constrained problems.  %In \cite{rinaldi2024stochastic} a trust-region method with convergence guarantees has been proposed. The method employs a tail-bound condition which is weaker than than employed in 
In \cite{dzahini2022expected} a stochastic direct search method is proposed for which both  asymptotic convergence and complexity analysis are carried out under standard probabilistic conditions.

\subsection{Our contribution}
We propose a derivative-free algorithm which is based on extrapolation techniques for the solution of problem \eqref{theoretical_prob}. For the proposed method, under quite standard assumptions, we manage to prove convergence to zero of the  norm of the gradient with probability one. Furthermore, we also give an iteration complexity result which aligns with those proposed in the recent literature, i.e. ${\cal O}(\eps^{-2})$.

\subsection{Assumptions}
In the sequel, we require the following assumptions.
\begin{assumption}\label{ass:compact}
The function $f(x)$ has compact level sets, i.e., for every $\alpha\in\Re$, the set
\[
    {\cal L}_\alpha = \{x\in\Re^n:\ f(x)\leq\alpha\}
\]
is compact. 
Furthermore, $f$ is bounded from above, i.e. $f_{\max}$ exists such that $f(x)\leq f_{\max}$ for all $x\in\Re^n$.
\end{assumption}

\begin{assumption}\label{ass:lipgrad}
The true objective function is continuously differentiable on $\Re^n$ with an L-Lipschitz continuous gradient, i.e. for all $x,y\in\Re^n$ it results
\[
\|\nabla f(x)-\nabla f(y)\| \leq L\|x-y\|.
\]
\end{assumption}

We note that, by Assumptions \ref{ass:compact}, the objective function $f$ is bounded from below on $\Re^n$, i.e. a value $f_{\rm low}$ exists such that $f_{\rm low}\leq f(x)$ for all $x\in\Re^n$.

Also, it is worth noticing that, even though our reference problem is (\ref{theoretical_prob}), we have only access to a black-box procedure that given a point $x\in\Re^n$, returns an observation of the real value $f(x)$, i.e. $F(x)$.

\section{The linesearch-type algorithm}
This section is devoted to the description of our algorithm, namely a Stochastic Derivative-free linesearch-based (SDFL) algorithm. The structure of the SDFL algorithm is inspired by the recent LAM algorithm studied in \cite{brilli2024worst}. In particular, the algorithm at every iteration carries out an exploration of the space around the current iterate $x_k$ using the coordinate direction. More specifically, the iteration starts by setting $y_k^1 = x_k$. Then, SDFL produces points $y_k^{i+1}$ for each $i=1,\dots,n$, such that $y_k^{i+1} = y_k^i + \alpha_k^ip_k^i$ with $\alpha_k^i \geq 0$. The stepsize $\alpha_k^i$ is equal to zero when a sufficient decrease along $\pm e_i$ cannot be obtained. On the other hand, $\alpha_k^i > 0$ is produced when a sufficient decrease can be obtained along either $e_i$ or $-e_i$. In more details, let $F(\cdot)$ denote  the estimate of the true objective function $f(\cdot)$ and $\bar\alpha_k^i$, $i=1,\dots,n$, denote the initial stepsizes for iteration $k$, then for every $i=1,\dots,n$, we have the following cases:
\begin{enumerate}
    \item if $F(y_k^i\pm\bar\alpha_k^ie_i) > F(y_k^i) -\gamma c\eps_f(\bar\alpha_k^i)^2$, then $\alpha_k^i = 0$;
    \item if $F(y_k^i+\bar\alpha_k^ip_k^i) \leq F(y_k^i) -\gamma c\eps_f(\bar\alpha_k^i)^2$, where $p_k^i = \pm e_i$, then $\alpha_k^i\geq\bar\alpha_k^i$ is computed such that either $\alpha_k^i = \bar\alpha_k^i$ and
    \[
    \begin{split}
    & F(y_k^i + \alpha_k^ip_k^i)  \leq F(y_k^i) - \gamma c\eps_f(\alpha_k^i)^2\\
    & F(y_k^i + 2\alpha_k^ip_k^i)  > F(y_k^i+\alpha_k^ip_k^i) - \gamma c\eps_f(\alpha_k^i)^2\\
    \end{split}
    \]
    or $\alpha_k^i >\bar\alpha_k^i$ and
    \[
    \begin{split}
    & F(y_k^i + \alpha_k^ip_k^i)  \leq F(y_k^i+\alpha_k^ip_k^i/2) - \gamma c\eps_f(\alpha_k^i/2)^2\\
    & F(y_k^i + 2\alpha_k^ip_k^i)  > F(y_k^i+\alpha_k^ip_k^i) - \gamma c\eps_f(\alpha_k^i)^2\\
    \end{split}
    \]
\end{enumerate}
When all the directions have been explored, so that points $y_k^i$, $i=1,\dots,n+1$, have been computed, the algorithm generates the tentative stepsizes for the next iteration, i.e. $\tilde\alpha_{k+1}^i$, $i=1,\dots,n$. In particular, if $x_{k+1} = x_k$, i.e. the iteration is deemed unsuccessful, 
\[
\tilde\alpha_{k+1}^i = \theta\bar\alpha_k^i,\quad i=1,\dots,n,
\]
otherwise, that is when $x_{k+1}\neq x_k$ and the iteration is successful, 
\[
\tilde\alpha_{k+1}^i = \max\{\alpha_k^i,\bar\alpha_k^i\},\quad i=1,\dots,n.
\]
 We report the scheme of the Stochastic derivative-free linesearch (SDFL) algorithm in the box below.
\par\medskip

\begin{algorithm}
    \caption{Stochastic Derivative-free Linesearch-based (SDFL)}
    \begin{algorithmic}[1]
        \State   {\bf data} $\theta, \in (0,1)$, $\gamma>2$, $c>0$, $\eta>0$, $\epsilon_f>0$, $x_0\in X$, $\tilde\alpha^i_0 > 0$, $i=1,\dots,n$, $\texttt{nF}=0$.

    \For{$k=0,1,\dots$}
    \State set $y_k^1 = x_k$
    \For{$i=1,2,\dots,n$}
    \State set \fbox{{\tt success = True}}, $\bar\alpha_k^i = \max\{\tilde\alpha_k^i,\eta\max_j\{\tilde\alpha_k^j\}\}$
    \State compute $F(y_k^i)$, $F(y_k^i+\bar\alpha_k^ie^i)$, set $\texttt{nF} = \texttt{nF}+2$
    \If{$F(y_k^i+\bar\alpha_k^ie^i)-F(y_k^i)> -\gamma c\eps_f (\bar\alpha_k^i)^2$}
    \State \texttt{/* try opposite direction */}
    \State compute  $F(y_k^i-\bar\alpha_k^ie^i)$, set $\texttt{nF} = \texttt{nF}+1$
    \If{$F(y_k^i-\bar\alpha_k^ie^i)-F(y_k^i)> -\gamma c\eps_f (\bar\alpha_k^i)^2$}
    \State \texttt{/* failure */}
    \State set $\alpha_k^i = 0, y_k^{i+1} = y_k^i$, \fbox{\tt success = False}
    \Else 
    \State set $p_k^i = -e^i$
    \EndIf
    \Else
    \State set $p_k^i = e^i$
    \EndIf
    \If{\fbox{\tt success}}
    \State \texttt{/* line search along $p_k^i$ */}
    \State set $\beta = 2\bar\alpha_k^i, \alpha=\bar\alpha_k^i$, compute $F(y_k^i+\beta p_k^i)$, set $\texttt{nF} = \texttt{nF}+1$
    \While{$F(y_k^i+\beta p_k^i)-F(y_k^i+\alpha p_k^i)\leq -\gamma c\eps_f (\beta-\alpha)^2$}
    \State set $\alpha = \beta$, $\beta = 2\alpha$
    \State compute $F(y_k^i+\beta p_k^i)$, set $\texttt{nF} = \texttt{nF}+1$
    \EndWhile
    \State set $\alpha_k^i = \alpha, y_k^{i+1} = y_k^i + \alpha_k^i p_k^i$
    \EndIf
    \EndFor
    \State set $x_{k+1}=y_k^{n+1}$
    \If{$x_{k+1}=x_k$}
    \State set $\tilde\alpha_{k+1}^i = \theta\bar\alpha_k^i$ for all $i$
    \Else
    \State set $\tilde\alpha_{k+1}^i = \max\{\alpha_k^i,\bar\alpha_k^i\}$ for all $i$
    \EndIf
    \EndFor
    
    \end{algorithmic}
\end{algorithm}

\section{Probabilistic estimates of $f$}

{We start this section by defining the accuracy of the function estimate $F(x)$ with respect to a precision parameter $\delta$. %(per Giampo: in realtà $c $ e $\eps_f$ compaiono sempre insieme non si capisce perchè li definiamo tutte e due. Uno si potrebbe togliere?)
}
\begin{defin}
Given $\eps_f>0$ and $c>0$, a point $x$ and a precision parameter $\delta$, we say that $F(x)$ is an $\eps_f$-accurate estimate of $f(x)$ when
\[
|F(x)-f(x)| \leq c\eps_f\delta^2.
\]
\end{defin}

The following result, which is proved (for instance) in \cite{dzahini2022expected}, relates sufficient decrease between estimates to sufficient decrease between true function values when estimates are $\eps_f$-accurate.

\begin{proposition}[{See \cite[Proposition 1]{dzahini2022expected}}]\label{reduce_epsaccurate}
Let $\gamma > 2$ and $F(x_i), F(x_{i+1})$ be $\epsilon_f$-accurate estimate of $f(x_i), f(x_{i+1})$, respectively. 

If $F(x_{i+1}) - F(x_i) \leq -\gamma c\eps_f\delta_i^2$, then
\[
f(x_{i+1}) - f(x_i) \leq -(\gamma-2)c\eps_f\delta_i^2.
\]
If $F(x_{i+1}) - F(x_i) > -\gamma c\eps_f\delta_i^2$, then
\[
f(x_{i+1}) - f(x_i) > -(\gamma+2)c\eps_f\delta_i^2.
\]

\end{proposition}

Considering algorithm SDFL, { in order to formalize the notion of conditioning on the past, we need to store all the information generated during the iterations  including the function evaluations produced in the line search. }
In particular, we define the following sets
\[
\begin{array}{l}
\phantom{\text{iteration $kk$\ }}
\begin{array}{rcl}{\cal G}_{-1,0} & = &\{F(\bm{x_0})\} %= {\cal G}_0 = {\cal G}_{0,0}
\end{array}\\
\text{iteration $0$}\left\{
\begin{array}{rcl}
{\cal G}_{-1,1} & = & {\cal G}_{-1,0}\cup\{F(\bm{x_{0,1}})\}\\
\vdots & &\qquad\vdots\\
{\cal G}_{-1,\ell_0} & = &{\cal G}_{-1,\ell_0-1}\cup \{F(\bm{x_{0,\ell_0}})\} = {\cal G}_{0} = {\cal G}_{0,0}
\end{array}\right.\\
\text{iteration $1$}\left\{
\begin{array}{rcl}
{\cal G}_{0,1} & = &{\cal G}_{0,0}\cup\{F(\bm{x_{1,1}})\}\\
\vdots & &\qquad\vdots\\
{\cal G}_{0,\ell_1} & = &{\cal G}_{0,\ell_1-1}\cup \{F(\bm{x_{1,\ell_1}})\} = {\cal G}_{1} = {\cal G}_{1,0}
\end{array}\right.\\
\phantom{\text{iteration $kk$\quad \quad\,\! }}\begin{array}{rcl}
\vdots & &\qquad\vdots
\end{array}\\
\text{iteration $k$}\left\{
\begin{array}{rcl}
{\cal G}_{k-1,1} & = &{\cal G}_{k-1,0}\cup\{F(\bm{x_{k,1}})\}\\
\vdots & &\qquad\vdots\\
{\cal G}_{k-1,\ell_k} & = &{\cal G}_{k-1,\ell_k-1}\cup \{F(\bm{x_{k,\ell_k}})\} = {\cal G}_{k} = {\cal G}_{k,0}
\end{array}\right.
\end{array}
\]
where $\bm{x_0}$ is the random variable associated with the initial point and $\bm{x_{k,i}}$ for $i=1,\dots,\ell_k$ are the random variables associated with the points generated by the algorithm during the $k$-th iteration.  Then, let ${\cal F}_{k,i}$ be the $\sigma$-field generated by ${\cal G}_{k,i}$.

We denote by ${\cal F}_{k}$ the $\sigma$-field generated by all the function estimates computed by the algorithm up to iteration $k$, i.e.
\[
{\cal F}_{k} = {\cal F}_{k-1,\ell_{k}}
\]
Note that, $\bm{x_k}$ is ${\cal F}_{k-1}$-measurable,  i.e. $\mathbb{E}[\bm{x_k}|{\cal F}_{k-1}] = \bm{x_k}$. 
Then, we can introduce the following fundamental assumption that will be used in the convergence analysis. 

% In particular, we assume that
% \[
% \mathbb{E}[F(x_i)|{\cal F}_{k-1}] = f(x_k)
% \]
% for any $x_i$ produced during the $k$th iteration of SDFL.

\begin{assumption}\label{ass:betaprob} For some $\beta\in(0,1)$, $\eps_f > 0$ and $c>0$, for any $k\geq 0$ and $i=1,\dots,\ell_k$,
\[
    \begin{split}
%& {\color{green}\mathbb{P}(J_k|{\cal F}_{k-1}) = \mathbb{E}(\mathbb{1}_{J_k}|{\cal F}_{k-1}) \geq \beta,}\\
& \mathbb{P}\left(|F(\bm{x}_i)-f(\bm{x}_i)|\leq c\eps_f\bm{\delta}_k^2\Big|{\cal F}_{k-1,i-1}\right) \geq \beta,\\
        &\mathbb{E}\left(|F(\bm{x}_i)-f(\bm{x}_i)|^2\Big| {\cal F}_{k-1,i-1}\right)\leq c^2\eps_f^2(1-\beta)\bm{\delta}_k^4   %,\\ %\|x_i-x_{i+1}\|^4\\
%&\mathbb{E}\left(|F(x_j)-f(x_j)|^2\Big| {\cal F}_{k-1}\right)\leq \eps_f^2(1-\beta)\delta_k^4\\ %\|x_i-x_{i+1}\|^4\\    
\end{split}
\]    
where $\bm{x}_i$ %and $x_j$ 
is the v.a. corresponding to the $i$-th point produced %either $x_k$ or 
%any point %pair of points 
in the $k$-th iteration of the algorithm and $\bm{\delta}_k = \min_{i=1,\dots,n}\{\bm{\bar\alpha}_k^i\}$. 
\end{assumption}
{Assumption \ref{ass:betaprob} basically requires that there is a prefixed probability ($\beta$) that the function values computed by the Algorithm at every iteration are $\eps_f$-accurate estimates of the true function values on the respective points.  It also provides a bound for the variance of the estimates. Furthermore note that the argument $\bm{x}_i$ is fixed by conditioning and  also that
$\mathbb{E}[{\bm{\delta}_k}|{\cal F}_{k-1}] = \bm{\delta}_k$.
}

Following the recent literature (see e.g. \cite{dzahini2022expected}), estimates can easily be computed as we briefly recall here for the sake of completeness.

So, given a point $x_i$ generated by SDFL at iteration $k$, let $\theta_h$, $h=1,\dots,p_j$ be $p_j$ independent realizations of the random variable $\theta$, then we define 
\[
F(x_i) = \frac{1}{p_j}\sum_{h=1}^{p_j}f(x_i,\theta^i_h),\qquad 
% F(x_j) = \frac{1}{p_j}\sum_{h=1}^{p_j}f(x_i,\theta^j_h).
\]
First, let us note that 
\[
\mathbb{V}[F(x_i)] \leq \frac{\mathbb{V}[f(x,\theta)]}{p_j}\leq \frac{V}{p_j}
\]
where $V>0$ is a constant such that the variance of $f(x,\theta)$ satisfies $\mathbb{V}[f(x,\theta)]\leq V < +\infty$ for all $x\in\Re^n$.

By the Markov inequality we can write
\[
\mathbb{P}\left(|F(\bm{x}_i)-f(\bm{x}_i)|^2> c^2\eps_f^2\bm{\delta}_k^4\Big|{\cal F}_{k-1,i-1}\right) 
\leq \frac{\mathbb{E}\left(|F(\bm{x}_i)-f(\bm{x}_i)|^2\Big| {\cal F}_{k-1,i-1}\right)}{c^2\eps_f^2\bm{\delta}_k^4}
%\leq \frac{\mathbb{V}[F(x_i)]}{c^2\eps_f^2\delta_k^4}
\leq \frac{V/p_j}{c^2\eps_f^2\bm{\delta}_k^4}
\]
thus, if we impose
\[
\frac{V}{p_jc^2\eps_f^2\bm{\delta}_k^4} \leq (1-\beta),
\]
% Then, it also results that $\mathbb{E}[F(x_i)|{\cal F}_{k-1}] = f(x_i)$. Furthermore, the value of $p_j$, i.e. the number of repeated black-box evaluations necessary to satisfy Assumption \ref{ass:betaprob} with reference to a given steplength $\delta_k$ is
which can be accomplished by choosing
\begin{equation}\label{eqassbeta2}
p_j \geq \dfrac{V}{c^2\eps_f^2(1-{\beta})\bm{\delta}_k^4},
\end{equation}
we have that
%Then, from Assumption \ref{ass:betaprob}, 
the following inequalities holds
\[
\frac{\mathbb{E}\left(|F(\bm{x}_i)-f(\bm{x}_i)|^2\Big| {\cal F}_{k-1,i-1}\right)}{c^2\eps_f^2\bm{\delta}_k^4}\leq (1-\beta)
\]
and
\begin{eqnarray}
\nonumber&&\mathbb{P}\left(|F(\bm{x}_i)-f(\bm{x}_i)|^2> c^2\eps_f^2\bm{\delta}_k^4\Big|{\cal F}_{k-1,i-1}\right) < 1-\beta,\\
\nonumber&&\mathbb{P}\left(|F(\bm{x}_i)-f(\bm{x}_i)|> c\eps_f\bm{\delta}_k^2\Big|{\cal F}_{k-1,i-1}\right) < 1-\beta,\\
\label{eqassbeta1}&&\mathbb{P}\left(|F(\bm{x}_i)-f(\bm{x}_i)|\leq c\eps_f\bm{\delta}_k^2\Big|{\cal F}_{k-1,i-1}\right) \geq \beta.
\end{eqnarray}
Hence, by \eqref{eqassbeta1} and \eqref{eqassbeta2}, we have that Assumption \ref{ass:betaprob} is satisfied.

\par\medskip

\begin{remark}\label{remark6}
At the start of iteration $k$, we can decide the number of repeated calls of the oracle such that, given $\beta\in(0,1)$, it results
\[
\mathbb{P}\left(|F(\bm{x}_k)-f(\bm{x}_k)|\leq c\eps_f\bm{\delta}_k^2 | {\cal F}_{k-1}\right) \geq \beta.
\]
Then, whichever is the sigma field ${\cal F}_{k-1,i-1}$, $i=1,\dots,\ell_k$, with same (or larger) number of repeated calls of the oracle, we can guarantee that
\[
\mathbb{P}\left(|F(\bm{x}_i)-f(\bm{x}_i)|\leq c\eps_f\bm{\delta}_k^2\Big|{\cal F}_{k-1,i-1}\right) \geq \beta.
\]
% Furthermore, we note that, for all $i=1,\dots,\ell_k$,
% \[
% \mathbb{P}\left(|F(x_i)-f(x_i)|\leq \eps_f\delta_k^2\Big|{\cal F}_{k}\right) = \mathbb{P}\left(|F(x_i)-f(x_i)|\leq \eps_f\delta_k^2\Big|{\cal F}_{k-1,i-1}\right) \geq \beta.
% \]
\end{remark}

\section{Convergence analysis of SDFL}

{In this section we consider the asymptotic properties of the sequence $\{x_k\}$ of points produced by the SDFL Algorithm.
\par\medskip\noindent
}

\subsection{Preliminary properties}
We first report an important result stating that the expansion step of Algorithm~SDFL is well-defined almost surely.

\begin{proposition}\label{prop:ls_welldefined}
Algorithm SDFL is well-defined, that is, the expansion step (when executed) always terminates in a finite number of steps almost surely.    
\end{proposition}

{\bf Proof}. 
{ 
We consider a generic $k$-th iteration of algorithm SDFL and we also consider the exploration performed by the algorithm starting from point $y_k^i$ along the direction $d_k^i$ which can either be $e_i$ or $-e_i$, for a generic index $i\in\{1,\dots,n\}$. Let us denote by $z_k^j$, $j=1,\dots$, the points generated by the expansion step (where $z_k^1=y_k^i$ is the initial point). Along with these points, the method computes the function values $F(z_k^1), F(z_k^2), \dots$ which build up ${\cal G}_{k-1,i}$ that, in turn, generates the events ${\cal F}_{k-1,i}$.
\par\medskip\noindent
To begin with, note that the points $\bm{z}_k^j$, for $j\geq 1$, are  
\[
\bm{z}_k^j = \bm{y}_k^i+ 2^{j-1}\bm{\bar\alpha}_k^i \bm{d}_k^i,\quad j=1,2,\dots.
\]
% Then, we have 
% \[
% \begin{split}
% F(x_1) & \leq F(x_0) -\gamma c\epsilon_f \bar\alpha^2 \\
% F(x_j) & \leq F(x_{j-1}) -\gamma c\epsilon_f 2^{2(j-2)}\bar\alpha^2,\qquad j=2,3,\dots,i.  
% \end{split}
% \]

Let us denote by ${\cal F}_{k-1,\ell_j}$ the $\sigma$-algebra generated by all the function values computed by the algorithm during iteration $k$ up to (and including) point $z_k^j$, and note that 
{
\begin{eqnarray*}
&&
\mathbb{P}(\{\bm{z}_k^{j+1}\ \text{is accepted}\}|{\cal F}_{k-1,\ell_j}) = \mathbb{P}\left(F(\bm{z}_k^{j+1}) - F(\bm{z}_k^{j}) \leq -\gamma c\epsilon_f 2^{2(j-1)}(\bm{\bar\alpha}_k^i)^2\big|{\cal F}_{k-1,\ell_j}\right)\\
&&= \mathbb{P}\Big(F(\bm{z}_k^{j}) - f(\bm{z}_k^{j})- F(\bm{z}_k^{j+1})+f(\bm{z}_k^{j+1}) \geq \gamma c\epsilon_f 2^{2(j-1)}(\bm{\bar\alpha}_k^i)^2-f(\bm{z}_k^{j})+f(\bm{z}_k^{j+1})\Big| {\cal F}_{k-1,\ell_j}\Big)\\
&&\le \mathbb{P}\Big(|F(\bm{z}_k^{j}) - f(\bm{z}_k^{j})- F(\bm{z}_k^{j+1})+f(\bm{z}_k^{j+1}) |\geq \gamma c\epsilon_f 2^{2(j-1)}(\bm{\bar\alpha}_k^i)^2-f(\bm{z}_k^{j})+f(\bm{z}_k^{j+1})\Big| {\cal F}_{k-1,\ell_j}\Big).
\end{eqnarray*}
By Assumption \ref{ass:compact} the true objective function $f$ is bounded both from below and above so we have that 
\begin{eqnarray*}
&&
\mathbb{P}(\{\bm{z}_k^{j+1}\ \text{is accepted}\}|{\cal F}_{k-1,\ell_j}) \le \\
&&\quad \mathbb{P}\Big(|F(\bm{z}_k^{j}) - f(\bm{z}_k^{j})- F(\bm{z}_k^{j+1})+f(\bm{z}_k^{j+1}) |\geq \gamma c\epsilon_f 2^{2(j-1)}(\bm{\bar\alpha}_k^i)^2-f_{\max}+f_{\rm low})\Big| {\cal F}_{k-1,\ell_j}\Big)
\end{eqnarray*}
Now for $j$ sufficiently large it is
$$ \gamma c\epsilon_f 2^{2(j-1)}(\bm{\bar\alpha}_k^i)^2-f_{\max}+f_{\rm low} >0$$
so that for the Chebycev inequality we have:
\begin{eqnarray*}
&&\mathbb{P}(\{\bm{z}_k^{j+1}\ \text{is accepted}\}|{\cal F}_{k-1,\ell_j}) \le \\  
&&\qquad \frac{\mathbb{E}\Big((F(\bm{z}_k^{j}) - f(\bm{z}_k^{j}))^2|{\cal F}_{k-1,\ell_j})\Big)+\mathbb{E}\Big((F(\bm{z}_k^{j+1}) - f(\bm{z}_k^{j+1}))^2|{\cal F}_{k-1,\ell_j})\Big)}{(\gamma c\epsilon_f 2^{2(j-1)}(\bm{\bar\alpha}_k^i)^2-f_{\max}+f_{\rm low})^2}\leq \\
&&\qquad \frac{2c^2\eps_f^2(1-\beta)\bm{\delta}_k^4}{(\gamma c\epsilon_f 2^{2(j-1)}(\bm{\bar\alpha}_k^i)^2-f_{\max}+f_{\rm low})^2},
\end{eqnarray*}
where the last inequality follows from Assumption \ref{ass:betaprob}.}}
Summing up for $j = 1,2,\dots$, we have
\[
\begin{split}
& \sum_{j=1}^\infty \mathbb{P}(\{\bm{z}_k^{j+1}\ \text{is accepted}\}|{\cal F}_{k-1,\ell_j}) \leq \sum_{j=1}^\infty \frac{2c^2\eps_f^2(1-\beta)\bm{\delta}_k^4}{(\gamma c\epsilon_f 2^{2(j-1)}(\bm{\bar\alpha}_k^i)^2-f_{\max}+f_{\rm low})^2}.%\\
%& \qquad \leq \sum_{i=1}^\infty \frac{f(x_i)-f(x_{i+1})}{\gamma c\epsilon_f 2^{2(i-1)}\bar\alpha^2}
\end{split}
\]
Then the series on the right hand side is convergent, thus the series 
\[
\sum_{j=1}^\infty\mathbb{P}(\{\bm{z}_k^{j+1}\ \text{is accepted}\}|{\cal F}_{k-1,\ell_j})
\]
is convergent too. Then, the proof is concluded by the Borel-Cantelli theorem, i.e. 
\[
\mathbb{P}\left(\Big\{\{\bm{z}_k^{j+1}\ \text{is accepted}\}|{\cal F}_{k-1,\ell_j}\Big\}\  i.o.\right) = 0, %\quad\text{for $i$ sufficiently large},
\]
meaning that, for $j$ sufficiently large, the expansion step produces a failure with probability one.
$\hfill\Box$
\par\medskip

\noindent
Now, with reference to iteration $k$, let us define the following event
\[
J_{k,i} = \{F(x_i)\ \text{is an $\eps_f$-accurate estimate of}\ f(x_i)\ \text{w.r.t.}\ \bm{\delta}_k\},
\]
and
\[
\mathbb{P}(J_{k,i}|{\cal F}_{k-1,i-1})\geq\beta.
\]
The following result is also proved in \cite{dzahini2022expected}. 
\begin{lemma}[See {\cite[Lemma 1]{dzahini2022expected}}]\label{lemma1_dzahini}
    Under Assumption \ref{ass:betaprob}, for all the points $\bm{x}_i$ produced in iteration $k$, the following holds true.
    \[
\mathbb{E}\left(\mathbb{1}_{\bar J_i}|F(\bm{x}_i)-f(\bm{x}_i)| \Big| {\cal F}_{k-1,i-1}\right) \leq c \epsilon_f (1-\beta)\bm{\delta}_k^2.
    \]
\end{lemma}

\noindent
{In order to characterize the evolution of the SDFL algorithm we introduce the following} improvement function:
\[
\Phi_k = \frac{\nu}{c\epsilon_f}(f(\bm{x}_k)-f_{\rm low}) + (1-\nu)\bm{\Delta}_k^2.
\]
where we denote by $\bm{\Delta}_k$ the maximum steplength at the $k$-th iteration, i.e.
\[
\bm{\Delta}_k = \max_{i=1,\dots,n}\{\bm{\tilde\alpha}_k^i\}.
\]
\begin{proposition}\label{phi_decrease}
Let $\nu$ be such that
\[
\frac{1}{1+(\gamma-2)(1/2)^2}< \nu<1,
\]
and $\beta\in(0,1)$ such that
\[
\frac{\beta^2}{1-\beta^2} > \frac{2\nu}{\min\{{\nu}(\gamma-2) \eta^2,(1-\nu)(1-\theta^2)\}}.
\]
The expected decrease in the random improvement function conditioned to the past satisfies, almost surely,
\begin{equation}\label{valbeta}
%\mathbb{E}[\Phi_{k+1} - \Phi_k | {\cal F}_{k-1,\ell_k-1}] \leq \Big[-\beta^2\min\{{\nu}(\gamma-2) \eta^2,(1-\nu)(1-\theta^2)\} + {\nu}(1-\beta^2)\Big]\bm{\Delta}_k^2
\mathbb{E}[\Phi_{k+1} - \Phi_k | {\cal F}_{k-1,\ell_k-1}] \leq -\frac{1}{2}\beta^2\min\{{\nu}(\gamma-2) \eta^2,(1-\nu)(1-\theta^2)\} \bm{\Delta}_k^2.
\end{equation}
%where the above expectation is negative provided that 
\end{proposition}

{\bf Proof.} Recall that, by Proposition~\ref{prop:ls_welldefined}, iteration $k$ is well-defined with probability one, i.e. iteration $k$ terminates with probability one. Then, under this situation, we start the proof by deriving some relations that hold true at successful iterations.
Hence, let us suppose that iteration $k$ is successful, i.e. $x_k\neq x_{k+1}$, and consider the following two cases
\begin{itemize}
    \item[(i)] $\Delta_k = \Delta_{k+1}$. Since the $k$-th iteration is of success, there is an index $\bar\jmath$ such that: 
	\[
	F(y_{k}^{\bar\jmath+1}) \leq F(y_{k}^{\bar\jmath}) - \gamma  c\epsilon_f(\bar\alpha_{k}^{\bar\jmath})^2
	\]
with
	\[
	\bar\alpha_{k}^{\bar\jmath} = \max\{\tilde\alpha_{k}^{\bar\jmath},\eta\Delta_{k}\} \geq \eta\Delta_{k},
	\]
	Then we have
	\[
	F(y_{k}^{\bar\jmath+1}) \leq F(y_{k}^{\bar\jmath}) - \gamma  c\epsilon_f(\bar\alpha_{k}^{\bar\jmath})^2 \leq F(y_{k}^{\bar\jmath}) - \gamma c\epsilon_f \eta^2\Delta_{k}^2.
	\]
	and, recalling that $F(y_k^{\bar\jmath})\leq F(x_k)$, we can write
	\begin{equation}\label{suffdeccond0}
	F(x_{k+1}) \leq F(x_k) - \gamma  c\epsilon_f \eta^2\Delta_{k}^2.
	\end{equation}

    \item[(ii)] $\Delta_{k+1} > \Delta_k$. We have that an index $\bar\jmath$ exists such that a linesearch has been performed along the $\bar\jmath$-th direction which determines a steplength  $\alpha_{k}^{\bar\jmath}$ satisfying
	\[
	\alpha_{k}^{\bar\jmath}=\tilde\alpha_{k+1}^{\bar\jmath}=\Delta_{k+1}.
	\]
	\par\noindent
 More specifically, we have
 \begin{equation}\label{suffdeccond1}
	\begin{split}
		F(x_{k+1}) &\leq F(y_{k}^{\bar\jmath} + \tilde\alpha_{k+1}^{\bar\jmath}d_{k+1}^{\bar\jmath}) \\
		&\leq F(y_{k}^{\bar\jmath}+\tilde\alpha_{k+1}^{\bar\jmath} d_{k+1}^{\bar\jmath}/2) - \gamma  c\epsilon_f(1/2)^2(\tilde\alpha_{k+1}^{\bar\jmath})^2 \\
        & \leq F(x_{k})- \gamma  c\epsilon_f(1/2)^2(\tilde\alpha_{k+1}^{\bar\jmath})^2\\
		&=F(x_{k})- \gamma  c\epsilon_f(1/2)^2\Delta_{k+1}^2.
	\end{split}
\end{equation}
\end{itemize}
% Hence, for successful iterations, we have
% \begin{equation}\label{decrFsuccess}
% F(x_k) \leq F(x_{k-1}) - \gamma c\epsilon_f \min\{(1-\delta)^2,\eta^2\}\Delta_k^2.
% \end{equation}

Now, let us consider the following event
\[
J_k = \{F(\bm{x}_k)\ \text{and}\ F(\bm{x}_{k+1})\ \text{are $\eps_f$ accurate estimates}\}.
\]
Note that $J_k$ is the conjunction of two  events whose  probability conditioned to ${\cal F}_{k-1,\ell_k-1}$ is at least $\beta$. Since the two events are independent conditionally to ${\cal F}_{k-1,\ell_k-1}$, recalling Assumption \ref{ass:betaprob} and {Remark \ref{remark6}}, we have 
\begin{equation}\label{prob_betasq}
\mathbb{P}(J_k|{\cal F}_{k-1,\ell_k-1})\geq \beta^2.
\end{equation}

We separately consider the cases of good and bad estimates.
\begin{enumerate}
    \item {\bf Good estimates} ($\mathbb{1}_{J_k}=1$), i.e. the $k$-th iteration is an iteration where estimates are good.
    \begin{itemize}
        \item[] \underline{Successful iteration}, $\bm{x}_{k+1}\neq \bm{x}_k$. 
        If $\bm{\Delta}_k = \bm{\Delta}_{k+1}$ ($\mathbb{1}_{S_k^=}=1$). Then, recalling \eqref{suffdeccond0} and   Proposition \ref{reduce_epsaccurate}, we have
        \[
        \mathbb{1}_{J_k}\mathbb{1}_{S_k^=}(f(\bm{x}_{k+1}) - f(\bm{x}_k))\leq  -\mathbb{1}_{J_k}\mathbb{1}_{S_k^=}(\gamma-2)c\epsilon_f \eta^2\bm{\Delta}_k^2.
        \]
        Then, 
        \begin{equation}\label{diffPhi_1}
        \mathbb{1}_{J_k}\mathbb{1}_{S_k^=}(\Phi_{k+1} -\Phi_k) = \mathbb{1}_{J_k}\mathbb{1}_{S_k^=}\frac{\nu}{c\epsilon_f}(f(\bm{x}_{k+1})-f(\bm{x}_k)) \leq -\mathbb{1}_{J_k}\mathbb{1}_{S_k^=}\nu(\gamma-2)\eta^2\bm{\Delta}_k^2.
        \end{equation}
        If, on the other hand, $\bm{\Delta}_{k+1}>\bm{\Delta}_k$ ($\mathbb{1}_{S_k^>}=1$), we have
        \begin{equation*}
        \begin{split}
                \mathbb{1}_{J_k}\mathbb{1}_{S_k^>}({\Phi}_{k+1} - {\Phi}_k) & = \mathbb{1}_{J_k}\mathbb{1}_{S_k^>}\frac{\nu}{c\epsilon_f}(f(\bm{x}_{k+1})-f(\bm{x}_k)) + (1-\nu)(\bm{\Delta}_{k+1}^2-\bm{\Delta}_k^2) %\\
                %& \leq \frac{\nu}{c\epsilon_f}(f(x_{k+1})-f(x_k)) + (1-\nu)\Delta_{k+1}^2
        \end{split}
        \end{equation*}
        Then, recalling \eqref{suffdeccond1} and Proposition \ref{reduce_epsaccurate}, we know that
        \[
        \mathbb{1}_{J_k}\mathbb{1}_{S_k^>}(f(\bm{x}_{k+1}) - f(\bm{x}_k)) \leq -\mathbb{1}_{J_k}\mathbb{1}_{S_k^>}(\gamma-2) c\epsilon_f (1/2)^2\bm{\Delta}_{k+1}^2
        \]
        so that
        \[
        \mathbb{1}_{J_k}\mathbb{1}_{S_k^>}(\Phi_{k+1} - \Phi_{k}) \leq \mathbb{1}_{J_k}\mathbb{1}_{S_k^>}(-{\nu}(\gamma-2) (1/2)^2\bm{\Delta}_{k+1}^2 + (1-\nu)\bm{\Delta}_{k+1}^2 -(1-\nu)\bm{\Delta}_k^2)
        \]
        Hence, when $\nu$ is sufficiently close to $1$, i.e. when $-{\nu}(\gamma-2) (1/2)^2 + (1-\nu) < 0$, that is
        \[
        \frac{1}{1+(\gamma-2)(1/2)^2}<\nu<1,
        \]
        we can write
        \begin{equation}\label{diffPhi_2}
        \mathbb{1}_{J_k}\mathbb{1}_{S_k^>}(\Phi_{k+1} - \Phi_k) \leq -\mathbb{1}_{J_k}\mathbb{1}_{S_k^>}(1-\nu)\bm{\Delta}_k^2.
        \end{equation}
        \item[] \underline{Unsuccessful iteration} ($\mathbb{1}_{\bar S_k}=1$), $x_{k+1}=x_k$ (and in this case $\bm{\Delta}_{k+1} = \theta\bm{\Delta}_k < \bm{\Delta}_k$). Then, we have
        \begin{equation}\label{diffPhi_3}
        \begin{split}
\mathbb{1}_{J_k}\mathbb{1}_{\bar S_k}(\Phi_{k+1} - \Phi_k) & = \mathbb{1}_{J_k}\mathbb{1}_{\bar S_k}(1-\nu)(\bm{\Delta}_{k+1}^2-\bm{\Delta}_k^2) \\
&= -\mathbb{1}_{J_k}\mathbb{1}_{\bar S_k}(1-\nu)(1-\theta^2)\bm{\Delta}_k^2.    
        \end{split}
        \end{equation}
        
    \end{itemize}
Then, in the case of good estimates, recalling \eqref{diffPhi_1}, \eqref{diffPhi_2}, and \eqref{diffPhi_3}, we can write
\begin{equation}\label{diffPhi_goodestimates}
\mathbb{1}_{J_k}(\Phi_{k+1} - \Phi_k) \leq -\mathbb{1}_{J_k}\min\{{\nu}(\gamma-2) \eta^2,(1-\nu)(1-\theta^2)\}\bm{\Delta}_k^2
\end{equation}
when $\nu$ is such that
\[
\frac{1}{1+(\gamma-2)(1/2)^2} < \nu < 1.
\]

Taking conditional expectation in \eqref{diffPhi_goodestimates} and recalling \eqref{prob_betasq}, we obtain
\begin{equation}\label{expect_diffPhi_good}
\mathbb{E}[\mathbb{1}_{J_k}(\Phi_{k+1} - \Phi_k) | {\cal F}_{k-1,\ell_k-1}] \leq -\beta^2\min\{{\nu}(\gamma-2) \eta^2,(1-\nu)(1-\theta^2)\}\bm{\Delta}_k^2
\end{equation}
    \item {\bf Bad estimates} ($\mathbb{1}_{\bar J_k}=1$)
    \begin{itemize}
        \item[] \underline{Successful iteration}, $\bm{x}_{k+1}\neq \bm{x}_k$ ($\mathbb{1}_{S_k}=1$). We can write
        \[
        \begin{split}
        &\mathbb{1}_{\bar J_k}\mathbb{1}_{S_k}(f(\bm{x}_{k+1}) - f(\bm{x}_k)) = \\
        & \quad \mathbb{1}_{\bar J_k}\mathbb{1}_{S_k}(F(\bm{x}_{k+1}) - F(\bm{x}_k) + f(\bm{x}_{k+1}) - F(\bm{x}_{k+1}) + F(\bm{x}_k) - f(\bm{x}_k)) \leq \\
        &  \quad \mathbb{1}_{\bar J_k}\mathbb{1}_{S_k}(F(\bm{x}_{k+1}) - F(\bm{x}_k) + |f(\bm{x}_{k+1}) - F(\bm{x}_{k+1})| + |F(\bm{x}_k) - f(\bm{x}_k)|) 
        \end{split}
        \]
        %Recalling %\eqref{suffdeccond1} when $\Delta_{k+1} > \Delta_k$, and \
        When $\bm{\Delta}_{k+1}=\bm{\Delta}_k$ ($\mathbb{1}_{S_k^=}=1$), recalling \eqref{suffdeccond0}, we can write 
        \[
        \begin{split}
        \mathbb{1}_{\bar J_k}\mathbb{1}_{S_k^=}(f(\bm{x}_{k+1}) - f(\bm{x}_k)) & \leq -\mathbb{1}_{\bar J_k}\mathbb{1}_{S_k^=}\gamma c\epsilon_f\eta^2\bm{\Delta}_k^2 + \\
        & \mathbb{1}_{\bar J_k}\mathbb{1}_{S_k^=}(|f(x_{k+1}) - F(x_{k+1})| + |f(x_k) - F(x_k)|)%\\
%        & \leq |f(x_k) - F(x_k)| + |f(x_{k-1}) - F(x_{k-1})|
        \end{split}
        \]
        Then, 
        \begin{equation}\label{diffPhi_4}
        \begin{split}
        \mathbb{1}_{\bar J_k}\mathbb{1}_{S_k^=}(\Phi_{k+1} - \Phi_k) &\leq \mathbb{1}_{\bar J_k}\mathbb{1}_{S_k^=}\frac{\nu}{c\epsilon_f}(|f(\bm{x}_{k+1}) - F(\bm{x}_{k+1})| + |f(\bm{x}_k) - F(\bm{x}_k)|) \\
        &-\mathbb{1}_{\bar J_k}\mathbb{1}_{S_k^=}\gamma \nu\eta^2\bm{\Delta}_k^2. 
        \end{split}
        \end{equation}
        When $\bm{\Delta}_{k+1} > \bm{\Delta}_k$ ($\mathbb{1}_{S_k^>}=1$), recalling \eqref{suffdeccond1}, we can write
        \[
        \begin{split}
        \mathbb{1}_{\bar J_k}\mathbb{1}_{S_k^>}(f(\bm{x}_{k+1}) - f(\bm{x}_k)) & \leq -\mathbb{1}_{\bar J_k}\mathbb{1}_{S_k^>}\gamma c\epsilon_f(1/2)^2\bm{\Delta}_{k+1}^2\\
        & + \mathbb{1}_{\bar J_k}\mathbb{1}_{S_k^>}(|f(\bm{x}_{k+1}) - F(\bm{x}_{k+1})| + |f(\bm{x}_k) - F(\bm{x}_k)|).%\\
        \end{split}
        \]
        Then,
        \[
        \begin{split}
            \mathbb{1}_{\bar J_k}\mathbb{1}_{S_k^>}(\Phi_{k+1} - \Phi_k) & \leq  \mathbb{1}_{\bar J_k}\mathbb{1}_{S_k^>}\frac{\nu}{c\epsilon_f}(|f(\bm{x}_{k+1}) - F(\bm{x}_{k+1})| + |f(\bm{x}_k) - F(\bm{x}_k)|) \\
            & - \mathbb{1}_{\bar J_k}\mathbb{1}_{S_k^>}(\gamma\nu (1/2)^2\bm{\Delta}_{k+1}^2 + (1-\nu)(\bm{\Delta}_{k+1}^2-\bm{\Delta}_k^2))
        \end{split}
        \]
        When $\nu$ is such that
        \[
        \frac{1}{1+\gamma(1/2)^2}< \nu<1
        \]
        we have
        \begin{equation}\label{diffPhi_5}
        \begin{split}
            \mathbb{1}_{\bar J_k}\mathbb{1}_{S_k^>}(\Phi_{k+1} - \Phi_k) & \leq  \mathbb{1}_{\bar J_k}\mathbb{1}_{S_k^>}\frac{\nu}{c\epsilon_f}(|f(\bm{x}_{k+1}) - F(\bm{x}_{k+1})| + |f(\bm{x}_k) - F(\bm{x}_k)|) \\
            & - \mathbb{1}_{\bar J_k}\mathbb{1}_{S_k^>}(1-\nu)\bm{\Delta}_k^2
        \end{split}
        \end{equation}
        
        \item[] \underline{Unsuccessful iteration}, $\bm{x}_{k+1}= \bm{x}_k$, $\bm{\Delta}_{k+1} = \theta\bm{\Delta}_k < \bm{\Delta}_k$ ($\mathbb{1}_{\bar S_k}$). We can write (as in the case of good estimates)        
        \begin{equation}\label{diffPhi_6}
        \begin{split}
        \mathbb{1}_{\bar J_k}\mathbb{1}_{\bar S_k}(\Phi_{k+1} - \Phi_k) & = \mathbb{1}_{\bar J_k}\mathbb{1}_{\bar S_k}(1-\nu)(\bm{\Delta}_{k+1}^2-\bm{\Delta}_k^2) = -\mathbb{1}_{\bar J_k}\mathbb{1}_{\bar S_k}(1-\nu)(1-\theta^2)\bm{\Delta}_k^2 \\
        &\leq\mathbb{1}_{\bar J_k}\mathbb{1}_{\bar S_k}\frac{\nu}{c\epsilon_f} (|f(\bm{x}_k) - F(\bm{x}_k)| + |f(\bm{x}_{k-1}) - F(\bm{x}_{k-1})|).
        \end{split}
        \end{equation}
    \end{itemize}    
    Then, in the case of bad estimates, recalling \eqref{diffPhi_4}, \eqref{diffPhi_5}, and \eqref{diffPhi_6}, we can always write
\begin{equation}\label{diffPhi_badestimates}
    \mathbb{1}_{\bar J_k}(\Phi_{k+1} - \Phi_k)\leq\mathbb{1}_{\bar J_k}\frac{\nu}{c\epsilon_f} (|f(\bm{x}_{k+1}) - F(\bm{x}_{k+1})| + |f(\bm{x}_k) - F(\bm{x}_k)|)
    \end{equation}
    when $\nu$ is such that
        \[
        \frac{1}{1+\gamma(1/2)^2}< \nu<1.
        \]
Taking conditional expectations in \eqref{diffPhi_badestimates} and recalling Lemma \ref{lemma1_dzahini}, this yields
\begin{equation}\label{expect_diffPhi_bad}
\mathbb{E}[\mathbb{1}_{\bar J_k}(\Phi_{k+1} - \Phi_k) | {\cal F}_{k-1,\ell_k-1}] \leq  2{\nu}(1-\beta^2)\bm{\delta}_k^2 \leq 2{\nu}(1-\beta^2)\bm{\Delta}_k^2.
\end{equation}
    
\end{enumerate}

Then, recalling \eqref{expect_diffPhi_good} and \eqref{expect_diffPhi_bad}, we can write
\[
\begin{split}
\mathbb{E}[\Phi_{k+1} - \Phi_k | {\cal F}_{k-1,\ell_k-1}] & = \mathbb{E}[(\mathbb{1}_{J_k} + \mathbb{1}_{\bar J_k})(\Phi_{k+1} - \Phi_k) | {\cal F}_{k-1,\ell_k-1}] \\
& \leq \Big[-\beta^2\min\{{\nu}(\gamma-2) \eta^2,(1-\nu)(1-\theta^2)\} + 2{\nu}(1-\beta^2)\Big]\bm{\Delta}_k^2\\
& \leq -\frac{1}{2}\beta^2\min\{{\nu}(\gamma-2) \eta^2,(1-\nu)(1-\theta^2)\} \bm{\Delta}_k^2
\end{split}
\]
where the second inequality follows from the requirement on $\beta$ that is
%above expectation is negative provided that 
\[
\frac{\beta^2}{1-\beta^2} > \frac{4\nu}{\min\{{\nu}(\gamma-2) \eta^2,(1-\nu)(1-\theta^2)\}}.
\]
Then, noting that
\[
\frac{1}{1+\gamma(1/2)^2}< \frac{1}{1+(\gamma-2)(1/2)^2},
\]
the proof is concluded.$\hfill\Box$
\par\medskip\noindent
By using the previous proposition, The following result shows the asymptotic properties of the sequence $\{\bm{\Delta}_k\}$.

\begin{theorem}\label{teo:limexpectdelta_tozero}Let $\beta$ be chosen according to (\ref{valbeta}) in Proposition \ref{phi_decrease}.
The sequence $\{\bm{\Delta}_k\}$ of maximum stepsizes produced by the algorithm is such that:
\begin{enumerate}
    \item[(i)] $\displaystyle\sum_{k=0}^\infty\mathbb{E}[\bm{\Delta}_{k}^2] < +\infty$;
    \item[(ii)]
    $\displaystyle\sum_{k=0}^\infty\bm{\Delta}_k^2 < +\infty$ {almost surely};
    \item[(iii)] $\displaystyle\lim_{k\to\infty}\mathbb{E}[\bm{\Delta}_k] = 0.$
\end{enumerate}
\end{theorem}
{\bf Proof}. From proposition \ref{phi_decrease}, we have that 
\[
\mathbb{E}[\Phi_{k+1}-\Phi_k|{\cal F}_{k-1,\ell_k-1}] \leq -\rho\bm{\Delta}_k^2,\quad\text{almost surely}
\]
where $\rho > 0$ (provided that $\beta^2$ sufficiently larger than $1/2$). Summing the above relation for $k=0,1,\dots,N$, we have
\[
\rho\sum_{k=0}^N\bm{\Delta}_k^2 \leq \sum_{k=0}^N\mathbb{E}[\Phi_{k}-\Phi_{k+1}|{\cal F}_{k-1,\ell_k-1}]. 
\]
Then, taking expectations on both sides and recalling that $\mathbb{E}[\mathbb{E}[\Phi_{k}-\Phi_{k+1}|{\cal F}_{k-1,\ell_k-1}] ] = \mathbb{E}[\Phi_k-\Phi_{k+1}]$ and that $\Phi$ is a non-negative function, we have
\[
\rho\sum_{k=0}^N\mathbb{E}[\bm{\Delta}_k^2] \leq \sum_{k=0}^N\mathbb{E}[\Phi_{k}-\Phi_{k+1}] = \mathbb{E}[\Phi_0] - \mathbb{E}[\Phi_{N+1}]\leq \mathbb{E}[\Phi_0]. 
\]
This implies
\[
\sum_{k=0}^\infty\mathbb{E}[\bm{\Delta}_k^2] \leq  \dfrac{\mathbb{E}[\Phi_0]}{\rho},
\]
which proves point (i). Then, we obtain
\begin{equation}\label{limexpectdelta2}
\lim_{k\to\infty}\mathbb{E}[\bm{\Delta}_k^2] = 0,
\end{equation}
and, reasoning as in \cite[Theorem 3]{dzahini2022expected}, 
\[
\sum_{k=0}^\infty\bm{\Delta}_k^2 < +\infty,\quad\text{almost surely},
\]
which proves point (ii).
Furthermore, recalling that
\[
\mathbb{E}[\bm{\Delta}_k^2] = \mathbb{V}[\bm{\Delta}_k] + \mathbb{E}[\bm{\Delta}_k]^2 \geq \mathbb{E}[\bm{\Delta}_k]^2,
\]
from \eqref{limexpectdelta2} we also have
\[
\lim_{k\to\infty}\mathbb{E}[\bm{\Delta}_k] = 0,
\]
finally proving point (iii) and concluding the proof. $\hfill\Box$

\par\medskip

% We also report a result from \cite[Theorem 2]{audet2021stochastic}.

% \begin{theorem}
% Let $\{x_k\}$ and $\{\Delta_k\}$ be the sequences of points and maximum stepsizes produced by the algorithm. Then, there exists at least an index set $K$ such that $\{x_k\}_K$ and $\{\Delta_k\}_K$ are almost surely convergent. In particular, it holds
% \[
% \begin{split}
%     & \lim_{k\to\infty,k\in K} x_k = \bar x, \\
%     & \lim_{k\to\infty,k\in K} \Delta_k = 0 
% \end{split}
% \]
% almost surely.
% \end{theorem}
\par\noindent

\subsection{Bound on gradient norm expectation and convergence}
{In the following proposition we bound the expected value of the square of the generic $i$th component of the gradient of $f$ at $x_k$ (i.e. $|\nabla f(x_k)^\top e_i|^2$) with the expected values of $\Delta_{k+1}^2$ and $\Delta_k^2$.}  %point out that the expected value of the square norm  of the gradient  $\nabla f(x_k)$ is bounded  by the expected value of $\Delta_{k+1}^2$

\begin{proposition}\label{bound_grad_prop}
Let Assumption \ref{ass:betaprob} hold. Then, for all $i=1,\dots,n$, we show that
\[
\mathbb{E}\left[|\nabla f(x_k)^\top e_i|^2\right] = {\cal O}(\mathbb{E}[\Delta_{k+1}^2]) + {\cal O}(\mathbb{E}[\Delta_{k}^2]). 
\]
\end{proposition}

% \begin{proposition}\label{bound_grad_prop}
% 	Suppose that Assumption \ref{ass:betaprob} holds. Let $\{x_k\}$ and $\{\Delta_k\}$ be the sequences produced by the SDFL algorithm. Then, the following bound holds
% \begin{eqnarray}
%     \mathbb{E}[\|\nabla f(x_k)\|] & \leq & \check c\mathbb{E}\left[\Delta_{k+1}\right] + 2\sqrt{n}\eps_f\sqrt{1-\beta}\mathbb{E}[\Delta_k],\label{bound_norm_grad}\\
%     \mathbb{E}[\|\nabla f(x_k)\|^2] & \leq &\check c_1\mathbb{E}[\Delta_{k+1}^2] + 4n\eps_f^2(1-\beta)\mathbb{E}[\Delta_{k}^2].\label{bound_norm1_grad} 
% \end{eqnarray}
% % 	Then, for each $k$ such that $x_{k+1}\neq  x_k$
% % 	\begin{equation}\label{bound1}
% % 	\|\nabla f(x_k)\| \leq \sqrt{n}\left(\frac{\eps_f(\gamma+2)+L(\sqrt{n}+1)}{{\delta}} \right)\max_{i=1,\dots,n}\{\tilde\alpha^i_{k+1}\},
% % \end{equation}
% % 	whereas, for each $k$ such that $x_{k+1} = x_k$
% % 	\begin{equation}\label{bound2}
% % 	\|\nabla f(x_k)\| \leq \sqrt{n}\frac{\eps_f(\gamma+2) + L}{\theta}\max_{i=1,\dots,n}\{\tilde\alpha^i_{k+1}\},
% % 	\end{equation}
% where the constants $\check c$ and $\check c_1$ only depend on the constants $\gamma,\theta$ defined in the SDFL algorithm, on the number of variables $n$ and on the Lipschitz constant $L$ of the gradient of the objective function.
% \end{proposition}

{\bf Proof.} 
% {\color{red}First of all, let us define the following event
% \[
% J_k^i = \{F(y_k^i+2\tilde\alpha_{k+1}^ie_i), F(y_k^i+\tilde\alpha_{k+1}^ie_i), \text{\ and\ } F(y_k^i+\tilde\alpha_{k+1}^ie_i/2)\ \text{are $\eps_f$ accurate} \}.
% \]
% Reasoning as before, 
% \[
% \mathbb{P}[J_k^i|{\cal F}_{k-1,\ell_k-1}] \geq \beta^3.
% \]}
We distinguish the following cases:
\begin{enumerate}
\item {\em Successful extrapolation} ($\mathbb{1}_{S_k^i}=1$), $\bm{\alpha}_k^i > 0$; then we have $\bm{\tilde\alpha}_{k+1}^i = \bm{\alpha}_k^i$ (we only consider the case in which direction $e_i$ is explored since similar reasonings apply to the case where the opposite direction $-e_i$ is explored). Then, we consider the event
\[
J_k^i = \{F(y_k^i+2\tilde\alpha_{k+1}^ie_i), F(y_k^i+\tilde\alpha_{k+1}^ie_i), \text{\ and\ } F(y_k^i+\tilde\alpha_{k+1}^ie_i/2)\ \text{are $\eps_f$ accurate} \}.
\]
and the following two subcases:
\begin{enumerate}
    \item good estimates ($\mathbb{1}_{J_k^i}=1$)
\[
\begin{split}
  & \mathbb{1}_{S_k^i}\mathbb{1}_{J_k^i}F\left(\bm{y}_{k}^i+{2\bm{\tilde\alpha}_{k+1}^i} e_i\right)> \mathbb{1}_{S_k^i}\mathbb{1}_{J_k^i}(F(\bm{y}_{k}^i+\bm{\tilde\alpha}_{k+1}^i e_i) - c\eps_f\gamma(\bm{\tilde\alpha}_{k+1}^i)^2),\\ %\quad
  & \mathbb{1}_{S_k^i}\mathbb{1}_{J_k^i}F(\bm{y}_{k}^i+ \bm{\tilde\alpha}_{k+1}^i e_i/2) \geq \mathbb{1}_{S_k^i}\mathbb{1}_{J_k^i}(F(\bm{y}_{k}^i + \bm{\tilde\alpha}_{k+1}^i e_i) + c\eps_f\gamma(1/2)^2(\bm{\tilde\alpha}_{k+1}^i)^2)
\end{split}
\]
Then, by Proposition \ref{reduce_epsaccurate}, we can write
\[
\begin{split}
  & \mathbb{1}_{S_k^i}\mathbb{1}_{J_k^i}f\left(\bm{y}_{k}^i+{2\bm{\tilde\alpha}_{k+1}^i} e_i\right)> \mathbb{1}_{S_k^i}\mathbb{1}_{J_k^i}(f(\bm{y}_{k}^i+\bm{\tilde\alpha}_{k+1}^i e_i) - c\eps_f(\gamma+2)(\bm{\tilde\alpha}_{k+1}^i)^2),\\ %\quad
  & \mathbb{1}_{S_k^i}\mathbb{1}_{J_k^i}f(\bm{y}_{k}^i+ \bm{\tilde\alpha}_{k+1}^i e_i/2) \geq \mathbb{1}_{S_k^i}\mathbb{1}_{J_k^i}(f(\bm{y}_{k}^i + \bm{\tilde\alpha}_{k+1}^i e_i) + c\eps_f(\gamma-2)(1/2)^2(\bm{\tilde\alpha}_{k+1}^i)^2).
\end{split}
\]
Then, from the Mean-Value Theorem we get,
\begin{eqnarray}\label{caseii_cond1a}
&& \phantom{-}\mathbb{1}_{S_k^i}\mathbb{1}_{J_k^i}\nabla f(\bm{\bar u}_k^i )^Te_i > -\mathbb{1}_{S_k^i}\mathbb{1}_{J_k^i}c\eps_f(\gamma+2)\bm{\tilde\alpha}_{k+1}^i,\\ \label{caseii_cond1b}
&&-\mathbb{1}_{S_k^i}\mathbb{1}_{J_k^i}\nabla f(\bm{\hat u}_k^i )^Te_i \geq \phantom{-}\mathbb{1}_{S_k^i}\mathbb{1}_{J_k^i}c\eps_f(\gamma-2)(1/2)\bm{\tilde\alpha}_{k+1}^i,
\end{eqnarray}
{where 
$\bm{\bar u}_k^i=\bm{y}_{k}^i+\bm{\bar\lambda}_k^i \bm{\tilde\alpha}_{k+1}^i e_i$ and
$\bm{\hat u}_k^i=y_{k}^i+\bm{\hat\lambda}_k^i (1/2)\bm{\tilde\alpha}_{k+1}^i e_i$, 
%$\bar v_k^i=y_{k}^i-\bar\mu_k^i\tilde\alpha_{k+1}^i e_i$, and
%$\hat v_k^i=y_{k}^i-\hat\mu_k^i(1/2)\tilde\alpha_{k+1}^i e_i$, 
with 
$\bm{\bar\lambda}_k^i, \bm{\hat\lambda}_k^i%, \bar\mu_k^i, \hat\mu_k^i 
\in (1,2)$.}

{}When (\ref{caseii_cond1a}) holds, we can write
\[
 \mathbb{1}_{S_k^i}\mathbb{1}_{J_k^i}[\nabla f(\bm{\bar u}_k^i ) - \nabla f(\bm{x}_k) + \nabla f(\bm{x}_k)]^Te_i > -\mathbb{1}_{S_k^i}\mathbb{1}_{J_k^i}c\eps_f(\gamma+2)\bm{\tilde\alpha}_{k+1}^i,
\]
so that we obtain
\begin{equation}\label{bound1_left}
\begin{split}
 \mathbb{1}_{S_k^i}\mathbb{1}_{J_k^i}\nabla f(\bm{x}_k)^Te_i & > -\mathbb{1}_{S_k^i}\mathbb{1}_{J_k^i}(c\eps_f(\gamma+2)\bm{\tilde\alpha}_{k+1}^i + L\|\bm{x}_k-\bm{\bar u}_k^i\| )\\
 & > -\mathbb{1}_{S_k^i}\mathbb{1}_{J_k^i}(c\eps_f(\gamma+2)\bm{\tilde\alpha}_{k+1}^i +L \|\bm{x}_k - \bm{y}_{k}^i\| + L \bm{\tilde\alpha}_{k+1}^i).
 \end{split}
\end{equation}
{}From (\ref{caseii_cond1b}), we can write
\[
 \mathbb{1}_{S_k^i}\mathbb{1}_{J_k^i}[\nabla f(\bm{\hat u}_k^i ) - \nabla f(\bm{x}_k) + \nabla f(\bm{x}_k)]^Te_i \leq -\mathbb{1}_{S_k^i}\mathbb{1}_{J_k^i} c\eps_f(\gamma-2)(1/2)\bm{\tilde\alpha}_{k+1}^i,
\]
so that, in this case, we obtain
\begin{equation}\label{bound1_right}
\begin{split}
 \mathbb{1}_{S_k^i}\mathbb{1}_{J_k^i}\nabla f(\bm{x}_k)^Te_i & \leq  - \mathbb{1}_{S_k^i}\mathbb{1}_{J_k^i}(c\eps_f(\gamma-2)(1/2)\bm{\tilde\alpha}_{k+1}^i - L\|\bm{x}_k-\bm{\hat u}_k^i\|) \\
 & \leq  \phantom{-}\mathbb{1}_{S_k^i}\mathbb{1}_{J_k^i} (c\eps_f(\gamma+2)\bm{\tilde\alpha}_{k+1}^i + L \|\bm{x}_k - \bm{y}_{k}^i\| + L \bm{\tilde\alpha}_{k+1}^i).
\end{split}
\end{equation}

Now, considering (\ref{bound1_left}) and (\ref{bound1_right}), we get
\begin{equation}\label{LAM_utile1}
 \mathbb{1}_{S_k^i}\mathbb{1}_{J_k^i}|\nabla f(\bm{x}_k)^T e_i| \leq \mathbb{1}_{S_k^i}\mathbb{1}_{J_k^i}\left({c\eps_f(\gamma+2)+L(\sqrt{n}+1)} \right)\bm{\Delta}_{k+1}.
\end{equation}

    \item bad estimates ($\mathbb{1}_{\bar J_k^i}=1$)
\[
\begin{split}
  & \mathbb{1}_{S_k^i}\mathbb{1}_{\bar J_k^i}F\left(\bm{y}_{k}^i+{2\bm{\tilde\alpha}_{k+1}^i} e_i\right)> \mathbb{1}_{S_k^i}\mathbb{1}_{\bar J_k^i}(F(\bm{y}_{k}^i+\bm{\tilde\alpha}_{k+1}^i e_i) - c\eps_f\gamma(\bm{\tilde\alpha}_{k+1}^i)^2),\\ %\quad
  & \mathbb{1}_{S_k^i}\mathbb{1}_{\bar J_k^i}F(\bm{y}_{k}^i+ \bm{\tilde\alpha}_{k+1}^i e_i/2) \geq \mathbb{1}_{S_k^i}\mathbb{1}_{\bar J_k^i}(F(\bm{y}_{k}^i + \bm{\tilde\alpha}_{k+1}^i e_i) + c\eps_f\gamma(1/2)^2(\bm{\tilde\alpha}_{k+1}^i)^2)
\end{split}
\]
Then, we can write
\[
\begin{split}
  & \mathbb{1}_{S_k^i}\mathbb{1}_{\bar J_k^i}(f\left(\bm{y}_{k}^i+{2\bm{\tilde\alpha}_{k+1}^i} e_i\right) + |F\left(\bm{y}_{k}^i+{2\bm{\tilde\alpha}_{k+1}^i} e_i\right)-f\left(\bm{y}_{k}^i+{2\bm{\tilde\alpha}_{k+1}^i} e_i\right)| ) \geq \\
  &\mathbb{1}_{S_k^i}\mathbb{1}_{\bar J_k^i}F\left(\bm{y}_{k}^i+{2\bm{\tilde\alpha}_{k+1}^i} e_i\right)> \mathbb{1}_{S_k^i}\mathbb{1}_{\bar J_k^i}(F(\bm{y}_{k}^i+\bm{\tilde\alpha}_{k+1}^i e_i) - c\eps_f\gamma(\bm{\tilde\alpha}_{k+1}^i)^2)\geq\\ %\quad
  & \mathbb{1}_{S_k^i}\mathbb{1}_{\bar J_k^i}(f(\bm{y}_{k}^i+\bm{\tilde\alpha}_{k+1}^i e_i)- |F(\bm{y}_{k}^i+\bm{\tilde\alpha}_{k+1}^i e_i)-f(\bm{y}_{k}^i+\bm{\tilde\alpha}_{k+1}^i e_i)| - c\eps_f\gamma(\bm{\tilde\alpha}_{k+1}^i)^2),
  \end{split}
  \]
  \[
  \begin{split}
  & \mathbb{1}_{S_k^i}\mathbb{1}_{\bar J_k^i}(f(\bm{y}_{k}^i+ \bm{\tilde\alpha}_{k+1}^i e_i/2) + |F(\bm{y}_{k}^i+ \bm{\tilde\alpha}_{k+1}^i e_i/2)-f(\bm{y}_{k}^i+ \bm{\tilde\alpha}_{k+1}^i e_i/2)|)\geq \\
  & \mathbb{1}_{S_k^i}\mathbb{1}_{\bar J_k^i}F(\bm{y}_{k}^i+ \bm{\tilde\alpha}_{k+1}^i e_i/2) \geq \mathbb{1}_{S_k^i}\mathbb{1}_{\bar J_k^i}(F(\bm{y}_{k}^i + \bm{\tilde\alpha}_{k+1}^i e_i) + c\eps_f\gamma(1/2)^2(\bm{\tilde\alpha}_{k+1}^i)^2)\geq\\
  & \mathbb{1}_{S_k^i}\mathbb{1}_{\bar J_k^i}(f(\bm{y}_{k}^i + \bm{\tilde\alpha}_{k+1}^i e_i) + c\eps_f\gamma(1/2)^2(\bm{\tilde\alpha}_{k+1}^i)^2 - |F(\bm{y}_{k}^i + \bm{\tilde\alpha}_{k+1}^i e_i)-f(\bm{y}_{k}^i + \bm{\tilde\alpha}_{k+1}^i e_i)|)
\end{split}
\]
Now, by introducing the following quantities,
\[
\begin{split}
    & \Delta F_{k,+}^i = |F\left(\bm{y}_{k}^i+{2\bm{\tilde\alpha}_{k+1}^i} e_i\right)-f\left(\bm{y}_{k}^i+{2\bm{\tilde\alpha}_{k+1}^i} e_i\right)|, \\
    & \Delta F_{k,-}^i = |F(\bm{y}_{k}^i+ \bm{\tilde\alpha}_{k+1}^i e_i/2)-f(\bm{y}_{k}^i+ \bm{\tilde\alpha}_{k+1}^i e_i/2)|,\\
    & \Delta F_{k,0}^i = |F(\bm{y}_{k}^i + \bm{\tilde\alpha}_{k+1}^i e_i)-f(\bm{y}_{k}^i + \bm{\tilde\alpha}_{k+1}^i e_i)|.
\end{split}
\]
the above relations can be rewritten as
\[
\begin{split}
  & \mathbb{1}_{S_k^i}\mathbb{1}_{\bar J_k^i}f\left(\bm{y}_{k}^i+{2\bm{\tilde\alpha}_{k+1}^i} e_i\right) \geq \\
  & \quad\quad \mathbb{1}_{S_k^i}\mathbb{1}_{\bar J_k^i}(f(\bm{y}_{k}^i+\bm{\tilde\alpha}_{k+1}^i e_i) - c\eps_f\gamma(\bm{\tilde\alpha}_{k+1}^i)^2- \Delta F_{k,0}^i - \Delta F_{k,+}^i),\\
  & \mathbb{1}_{S_k^i}\mathbb{1}_{\bar J_k^i}f(\bm{y}_{k}^i+ \bm{\tilde\alpha}_{k+1}^i e_i/2) \geq \\
  & \quad\quad \mathbb{1}_{S_k^i}\mathbb{1}_{\bar J_k^i}(f(\bm{y}_{k}^i + \bm{\tilde\alpha}_{k+1}^i e_i) + c\eps_f\gamma(1/2)^2(\bm{\tilde\alpha}_{k+1}^i)^2 - \Delta F_{k,0}^i-\Delta F_{k,-}^i)
\end{split}
\]
Then, from the Mean-Value Theorem we get,
\begin{eqnarray}\label{caseiii_cond1a}
\phantom{-}\mathbb{1}_{S_k^i}\mathbb{1}_{J_k^i}\nabla f(\bm{\bar u}_k^i )^Te_i &\geq& -\mathbb{1}_{S_k^i}\mathbb{1}_{J_k^i}(c\eps_f\gamma\bm{\tilde\alpha}_{k+1}^i + \frac{\Delta F_{k,0}^i + \Delta F_{k,+}^i}{\bm{\tilde\alpha}_{k+1}^i}) \\ \nonumber &\geq& -\mathbb{1}_{S_k^i}\mathbb{1}_{J_k^i}(c\eps_f\gamma\bm{\tilde\alpha}_{k+1}^i + \frac{\Delta F_{k,0}^i + \Delta F_{k,+}^i+\Delta F_{k,-}^i}{\bm{\tilde\alpha}_{k+1}^i})\\ \label{caseiii_cond1b}
-\mathbb{1}_{S_k^i}\mathbb{1}_{J_k^i}\nabla f(\bm{\hat u}_k^i )^Te_i &\geq& \phantom{-}\mathbb{1}_{S_k^i}\mathbb{1}_{J_k^i}(c\eps_f\gamma(1/2)\bm{\tilde\alpha}_{k+1}^i- \frac{\Delta F_{k,0}^i+\Delta F_{k,-}^i}{\bm{\tilde\alpha}_{k+1}^i})\\
&\geq& -\mathbb{1}_{S_k^i}\mathbb{1}_{J_k^i}\nonumber(c\eps_f\gamma\bm{\tilde\alpha}_{k+1}^i + \frac{\Delta F_{k,0}^i + \Delta F_{k,+}^i+ \Delta F_{k,-}^i}{\bm{\tilde\alpha}_{k+1}^i})
\end{eqnarray}
{where 
$\bm{\bar u}_k^i=\bm{y}_{k}^i+\bm{\bar\lambda}_k^i \bm{\tilde\alpha}_{k+1}^i e_i$ and
$\bm{\hat u}_k^i=y_{k}^i+\bm{\hat\lambda}_k^i (1/2)\bm{\tilde\alpha}_{k+1}^i e_i$, 
%$\bar v_k^i=y_{k}^i-\bar\mu_k^i\tilde\alpha_{k+1}^i e_i$, and
%$\hat v_k^i=y_{k}^i-\hat\mu_k^i(1/2)\tilde\alpha_{k+1}^i e_i$, 
with 
$\bm{\bar\lambda}_k^i, \bm{\hat\lambda}_k^i%, \bar\mu_k^i, \hat\mu_k^i 
\in (1,2)$.}
Then, reasoning as in the previous case, considering that $\bm{\bar\alpha}_k^i\leq\bm{\tilde\alpha}_{k+1}^i$ and recalling that  $\bm{\delta}_k = \min_{i=1,\dots,n}\bm{\bar\alpha}_k^i$, we finally obtain
\begin{equation}\label{LAM_utile1i}
\begin{split}
 \mathbb{1}_{S_k^i}\mathbb{1}_{\bar J_k^i}|\nabla f(\bm{x}_k)^T e_i| & \leq \mathbb{1}_{S_k^i}\mathbb{1}_{\bar J_k^i}\left({c\eps_f\gamma+L(\sqrt{n}+1)} \right)\bm{\Delta}_{k+1}  \\ 
 &\qquad +\mathbb{1}_{S_k^i}\mathbb{1}_{\bar J_k^i}\frac{\Delta F_{k,0}^i + \Delta F_{k,+}^i+ \Delta F_{k,-}^i}{\bm{\delta}_{k}}.
\end{split}
\end{equation}

\end{enumerate}
\item {\em Unsuccessful extrapolation}, $\bm{\alpha}_k^i = 0$ ($\mathbb{1}_{\bar S_k^i}=1$). Then, we consider the event
\[
I_k^i = \{F(y_k^i-\bar\alpha_{k}^ie_i), F(y_k^i), \text{\ and\ } F(y_k^i+\bar\alpha_{k}^ie_i)\ \text{are $\eps_f$ accurate} \}.
\]
and the following two subcases:
\begin{enumerate}
    \item good estimates ($\mathbb{1}_{I_k^i}=1$); we have:
\begin{eqnarray*}
 && \mathbb{1}_{\bar S_k^i}\mathbb{1}_{I_k^i}F(\bm{y}_{k}^i+\bm{\bar\alpha}_k^i e_i)> \mathbb{1}_{\bar S_k^i}\mathbb{1}_{I_k^i}(F(\bm{y}_{k}^i)-c\eps_f\gamma(\bm{\bar\alpha}_k^i)^2),\\
 && \mathbb{1}_{\bar S_k^i}\mathbb{1}_{I_k^i}F(\bm{y}_{k}^i-\bm{\bar\alpha}_k^i e_i)> \mathbb{1}_{\bar S_k^i}\mathbb{1}_{I_k^i}(F(\bm{y}_{k}^i)-c\eps_f\gamma(\bm{\bar\alpha}_k^i)^2).
\end{eqnarray*}
which, recalling Proposition \ref{reduce_epsaccurate}, gives
\begin{eqnarray*}
 && \mathbb{1}_{\bar S_k^i}\mathbb{1}_{I_k^i}f(\bm{y}_{k}^i+\bm{\bar\alpha}_k^i e_i)> \mathbb{1}_{\bar S_k^i}\mathbb{1}_{I_k^i}(f(\bm{y}_{k}^i)-c\eps_f(\gamma+2)(\bm{\bar\alpha}_k^i)^2),\\
 && \mathbb{1}_{\bar S_k^i}\mathbb{1}_{I_k^i}f(\bm{y}_{k}^i-\bm{\bar\alpha}_k^i e_i)> \mathbb{1}_{\bar S_k^i}\mathbb{1}_{I_k^i}(f(\bm{y}_{k}^i)-c\eps_f(\gamma+2)(\bm{\bar\alpha}_k^i)^2).
\end{eqnarray*}
Then we get from the Mean-Value Theorem
\begin{eqnarray}
&& \mathbb{1}_{\bar S_k^i}\mathbb{1}_{I_k^i}\nabla f(\bm{u}_k^i )^Te_i > -\mathbb{1}_{\bar S_k^i}\mathbb{1}_{I_k^i}c\eps_f(\gamma+2)\bm{\bar\alpha}_k^i,\label{eq1}\\
&& \mathbb{1}_{\bar S_k^i}\mathbb{1}_{I_k^i}\nabla f(\bm{v}_k^i )^Te_i <  \phantom{-}\mathbb{1}_{\bar S_k^i}\mathbb{1}_{I_k^i}c\eps_f(\gamma+2)\bm{\bar\alpha}_k^i,\label{eq2}
\end{eqnarray}
where $\bm{u}_k^i=\bm{y}_{k}^i+\bm{\lambda}_k^i \bm{\bar\alpha}_k^i e_i$ and $\bm{v}_k^i=\bm{y}_{k}^i-\bm{\mu}_k^i\bm{\bar\alpha}_k^i e_i$ with $\bm{\lambda}_k^i, \bm{\mu}_k^i \in (0,1)$.
{}From (\ref{eq1}) and (\ref{eq2}) and the Lipschitz continuity of $\nabla f$, we have that
\begin{eqnarray*}
 \mathbb{1}_{\bar S_k^i}\mathbb{1}_{I_k^i}\nabla f(\bm{x}_k)^T e_i &>&          - \mathbb{1}_{\bar S_k^i}\mathbb{1}_{I_k^i}(c\eps_f(\gamma+2)\bm{\bar\alpha}_k^i +L \|\bm{x}_k - \bm{u}_k^i\|) \\
 & >&           - \mathbb{1}_{\bar S_k^i}\mathbb{1}_{I_k^i}(c\eps_f(\gamma+2)\bm{\bar\alpha}_k^i +L \|\bm{x}_k - \bm{y}_{k}^i\| - L \bm{\bar\alpha}_k^i),\\
 \mathbb{1}_{\bar S_k^i}\mathbb{1}_{I_k^i}\nabla f(\bm{x}_k)^T e_i &<& \phantom{-}\mathbb{1}_{\bar S_k^i}\mathbb{1}_{I_k^i}(c\eps_f(\gamma+2)\bm{\bar\alpha}_k^i +L \|\bm{x}_k - \bm{v}_k^i\|) \\
 &<& \phantom{-}\mathbb{1}_{\bar S_k^i}\mathbb{1}_{I_k^i}(c\eps_f(\gamma+2)\bm{\bar\alpha}_k^i +L \|\bm{x}_k - \bm{y}_{k}^i\| + L \bm{\bar\alpha}_k^i).
\end{eqnarray*}
Recalling that, from the instructions of the algorithm, we have $\bm{\bar\alpha}_k^i \leq \bm{\tilde\alpha}_{k+1}^i/\theta$, we can write
\[
\begin{split}
 & \mathbb{1}_{\bar S_k^i}\mathbb{1}_{I_k^i}|\nabla f(\bm{x}_k)^T e_i| < \mathbb{1}_{\bar S_k^i}\mathbb{1}_{I_k^i}(c\eps_f(\gamma+2)+L)\bm{\bar\alpha}_k^i +\mathbb{1}_{\bar S_k^i}\mathbb{1}_{I_k^i}L \|\bm{x}_k - \bm{y}_{k}^i\| \\
 &\qquad \leq \mathbb{1}_{\bar S_k^i}\mathbb{1}_{I_k^i}(c\eps_f(\gamma+2)+L)\bm{\bar\alpha}_k^i +\mathbb{1}_{\bar S_k^i}\mathbb{1}_{I_k^i}L\sqrt{n}\max_{i=1,\dots,n}\{\bm{\tilde\alpha}_{{k+1}}^i\}  \\ 
 &\qquad \leq \mathbb{1}_{\bar S_k^i}\mathbb{1}_{I_k^i}(c\eps_f(\gamma+2)+L)\ {\bm{\tilde\alpha}_{k+1}^i} +\mathbb{1}_{\bar S_k^i}\mathbb{1}_{I_k^i}L\sqrt{n}\max_{i=1,\dots,n}\{\bm{\tilde\alpha}_{{k+1}}^i\}  \\ 
 &\qquad \leq\mathbb{1}_{\bar S_k^i}\mathbb{1}_{I_k^i}{(c\eps_f(\gamma+2)+L)}\max_{i=1,\dots,n}\{\bm{\tilde\alpha}_{k+1}^i\} +\mathbb{1}_{\bar S_k^i}\mathbb{1}_{I_k^i}L\sqrt{n}\max_{i=1,\dots,n}\{\bm{\tilde\alpha}_{{k+1}}^i\},
 \end{split}
\]
so that
%and, since $\max_{i=1,\dots,n}\{\bar\alpha_k^i\} = \max_{i=1,\dots,n}\{\tilde\alpha_k^i\}$,
\begin{equation}\label{case_b_nuovo}
 \mathbb{1}_{\bar S_k^i}\mathbb{1}_{I_k^i}|\nabla f(\bm{x}_k)^T e_i| \leq \mathbb{1}_{\bar S_k^i}\mathbb{1}_{I_k^i}\biggl({c\eps_f(\gamma+2)+L(\sqrt{n}+1)} \biggr)\bm{\Delta}_{k+1}.
\end{equation}

    \item bad estimates ($\mathbb{1}_{\bar I_k^i}=1$)

\begin{eqnarray*}
 && \mathbb{1}_{\bar S_k^i}\mathbb{1}_{\bar I_k^i}F(\bm{y}_{k}^i+\bm{\bar\alpha}_k^i e_i)> \mathbb{1}_{\bar S_k^i}\mathbb{1}_{\bar I_k^i}(F(\bm{y}_{k}^i)-c\eps_f\gamma(\bm{\bar\alpha}_k^i)^2),\\
 && \mathbb{1}_{\bar S_k^i}\mathbb{1}_{\bar I_k^i}F(\bm{y}_{k}^i-\bm{\bar\alpha}_k^i e_i)> \mathbb{1}_{\bar S_k^i}\mathbb{1}_{\bar I_k^i}(F(\bm{y}_{k}^i)-c\eps_f\gamma(\bm{\bar\alpha}_k^i)^2).
\end{eqnarray*}
and we can write
\begin{eqnarray*}
 &&   \mathbb{1}_{\bar S_k^i}\mathbb{1}_{\bar I_k^i}(f(\bm{y}_{k}^i+\bm{\bar\alpha}_k^i e_i) + |F(\bm{y}_{k}^i+\bm{\bar\alpha}_k^i e_i)-f(\bm{y}_{k}^i+\bm{\bar\alpha}_k^i e_i)|)> \mathbb{1}_{\bar S_k^i}\mathbb{1}_{\bar I_k^i}F(\bm{y}_{k}^i+\bm{\bar\alpha}_k^i e_i)  > \\
 && \qquad \mathbb{1}_{\bar S_k^i}\mathbb{1}_{\bar I_k^i}(F(\bm{y}_{k}^i)-c\eps_f\gamma(\bm{\bar\alpha}_k^i)^2) > \mathbb{1}_{\bar S_k^i}\mathbb{1}_{\bar I_k^i}(f(\bm{y}_{k}^i) - |F(\bm{y}_{k}^i) - f(\bm{y}_{k}^i)| -c\eps_f\gamma(\bm{\bar\alpha}_k^i)^2),\\
 && \mathbb{1}_{\bar S_k^i}\mathbb{1}_{\bar I_k^i}(f(\bm{y}_{k}^i-\bm{\bar\alpha}_k^i e_i) + |F(\bm{y}_{k}^i-\bm{\bar\alpha}_k^i e_i)-f(\bm{y}_{k}^i-\bm{\bar\alpha}_k^i e_i)|)> \mathbb{1}_{\bar S_k^i}\mathbb{1}_{\bar I_k^i}F(\bm{y}_{k}^i-\bm{\bar\alpha}_k^i e_i)  > \\
 && \qquad \mathbb{1}_{\bar S_k^i}\mathbb{1}_{\bar I_k^i}(F(\bm{y}_{k}^i)-c\eps_f\gamma(\bm{\bar\alpha}_k^i)^2) > \mathbb{1}_{\bar S_k^i}\mathbb{1}_{\bar I_k^i}(f(\bm{y}_{k}^i) - |F(\bm{y}_{k}^i) - f(\bm{y}_{k}^i)| -c\eps_f\gamma(\bm{\bar\alpha}_k^i)^2).
\end{eqnarray*}
Now, by introducing the following quantities,
\[
\begin{split}
    & \Delta G_{k,+}^i = |F\left(\bm{y}_{k}^i+{\bm{\bar\alpha}_{k}^i} e_i\right)-f\left(\bm{y}_{k}^i+{\bm{\bar\alpha}_{k}^i} e_i\right)|, \\
    & \Delta G_{k,-}^i = |F(\bm{y}_{k}^i- \bm{\bar\alpha}_{k}^i e_i)-f(\bm{y}_{k}^i- \bm{\bar\alpha}_{k}^i e_i)|,\\
    & \Delta G_{k,0}^i = |F(\bm{y}_{k}^i) -f(\bm{y}_{k}^i)|,
\end{split}
\]
we have
\begin{eqnarray*}
 && \mathbb{1}_{\bar S_k^i}\mathbb{1}_{\bar I_k^i}f(\bm{y}_{k}^i+\bm{\bar\alpha}_k^i e_i)> \mathbb{1}_{\bar S_k^i}\mathbb{1}_{\bar I_k^i}(f(\bm{y}_{k}^i)-c\eps_f\gamma(\bm{\bar\alpha}_k^i)^2 -\Delta G_{k,0}^i - \Delta G_{k,+}^i),\\
 && \mathbb{1}_{\bar S_k^i}\mathbb{1}_{\bar I_k^i}f(\bm{y}_{k}^i-\bm{\bar\alpha}_k^i e_i)> \mathbb{1}_{\bar S_k^i}\mathbb{1}_{\bar I_k^i}(f(\bm{y}_{k}^i)-c\eps_f\gamma(\bm{\bar\alpha}_k^i)^2 -\Delta G_{k,0}^i - \Delta G_{k,-}^i).
\end{eqnarray*}
Then we get from the Mean Value Theorem
\begin{eqnarray}
&& \mathbb{1}_{\bar S_k^i}\mathbb{1}_{\bar I_k^i}\nabla f(\bm{u}_k^i )^Te_i > -\mathbb{1}_{\bar S_k^i}\mathbb{1}_{\bar I_k^i}\left(c\eps_f\gamma\bm{\bar\alpha}_k^i\label{eq1} + \dfrac{\Delta G_{k,0}^i}{\bm{\bar\alpha}_k^i} + \dfrac{\Delta G_{k,+}^i}{\bm{\bar\alpha}_k^i}\right) \\
&& \mathbb{1}_{\bar S_k^i}\mathbb{1}_{\bar I_k^i}\nabla f(\bm{v}_k^i )^Te_i < \phantom{-}\mathbb{1}_{\bar S_k^i}\mathbb{1}_{\bar I_k^i}\left(c\eps_f\gamma\bm{\bar\alpha}_k^i\label{eq2} + \dfrac{\Delta G_{k,0}^i}{\bm{\bar\alpha}_k^i} + \dfrac{\Delta G_{k,-}^i}{\bm{\bar\alpha}_k^i}\right),
\end{eqnarray}
where $\bm{u}_k^i=\bm{y}_{k}^i+\bm{\lambda}_k^i \bm{\bar\alpha}_k^i e_i$ and $\bm{v}_k^i=\bm{y}_{k}^i-\bm{\mu}_k^i\bm{\bar\alpha}_k^i e_i$ with $\bm{\lambda}_k^i, \bm{\mu}_k^i \in (0,1)$. {}From (\ref{eq1}) and (\ref{eq2}) and the Lipschitz continuity of $\nabla f$, we have that
\begin{eqnarray*}
 && \mathbb{1}_{\bar S_k^i}\mathbb{1}_{\bar I_k^i}\nabla f(\bm{x}_k)^T e_i >  -\mathbb{1}_{\bar S_k^i}\mathbb{1}_{\bar I_k^i}\left( c\eps_f\gamma\bm{\bar\alpha}_k^i +L \|\bm{x}_k - \bm{y}_{k}^i\| + L \bm{\bar\alpha}_k^i + \dfrac{\Delta G_{k,0}^i}{\bm{\bar\alpha}_k^i} + \dfrac{\Delta G_{k,+}^i}{\bm{\bar\alpha}_k^i}\right),\\
 && \mathbb{1}_{\bar S_k^i}\mathbb{1}_{\bar I_k^i}\nabla f(\bm{x}_k)^T e_i <  \mathbb{1}_{\bar S_k^i}\mathbb{1}_{\bar I_k^i}\left( c\eps_f\gamma\bm{\bar\alpha}_k^i +L \|\bm{x}_k - \bm{y}_{k}^i\| + L \bm{\bar\alpha}_k^i + \dfrac{\Delta G_{k,0}^i}{\bm{\bar\alpha}_k^i} + \dfrac{\Delta G_{k,-}^i}{\bm{\bar\alpha}_k^i}\right).
\end{eqnarray*}
Then, recalling that $\bm{\delta}_k  = \min_{i=1,\dots,n}\bm{\bar\alpha}_k^i$, we have
\begin{eqnarray*}
 && \mathbb{1}_{\bar S_k^i}\mathbb{1}_{\bar I_k^i}\nabla f(\bm{x}_k)^T e_i >  -\mathbb{1}_{\bar S_k^i}\mathbb{1}_{\bar I_k^i}\left( c\eps_f\gamma\bm{\bar\alpha}_k^i +L \|\bm{x}_k - \bm{y}_{k}^i\| + L \bm{\bar\alpha}_k^i + \dfrac{\Delta G_{k,0}^i}{\bm{\delta}_k} + \dfrac{\Delta G_{k,+}^i}{\bm{\delta}_k}\right),\\
 && \mathbb{1}_{\bar S_k^i}\mathbb{1}_{\bar I_k^i}\nabla f(\bm{x}_k)^T e_i <  \mathbb{1}_{\bar S_k^i}\mathbb{1}_{\bar I_k^i}\left( c\eps_f\gamma\bm{\bar\alpha}_k^i +L \|\bm{x}_k - \bm{y}_{k}^i\| + L \bm{\bar\alpha}_k^i + \dfrac{\Delta G_{k,0}^i}{\bm{\delta}_k} + \dfrac{\Delta G_{k,-}^i}{\bm{\delta}_k}\right).
\end{eqnarray*}
Hence, we can write
\begin{eqnarray}\label{bound_grad_bad_casei}
\mathbb{1}_{\bar S_k^i}\mathbb{1}_{\bar I_k^i}|\nabla f(\bm{x}_k)^T e_i | &<& \mathbb{1}_{\bar S_k^i}\mathbb{1}_{\bar I_k^i}\Bigg((c\eps_f\gamma +L(\sqrt{n}  + 1)) \bm{\Delta}_{k+1} \\
&& \nonumber + \dfrac{\Delta G_{k,0}^i}{\bm{\delta}_k} + \dfrac{\Delta G_{k,+}^i}{\bm{\delta}_k}+\dfrac{\Delta G_{k,-}^i}{\bm{\delta}_k}\Bigg).
\end{eqnarray}
    
\end{enumerate}

\end{enumerate}

Now, let us denote 
\begin{equation}\label{def:chat}
\hat c \triangleq c\eps_f(\gamma+2)+L(\sqrt{n}+1).
\end{equation}
Hence, when $\mathbb{1}_{S_k^i}=1$, using \eqref{LAM_utile1} and \eqref{LAM_utile1i},
\begin{eqnarray*}
 \mathbb{1}_{S_k^i}\mathbb{1}_{J_k^i}|\nabla f(\bm{x}_k)^T e_i| &\leq& \mathbb{1}_{S_k^i}\mathbb{1}_{J_k^i}\hat c\bm{\Delta}_{k+1},\\
 \mathbb{1}_{S_k^i}\mathbb{1}_{\bar J_k^i}|\nabla f(\bm{x}_k)^T e_i| & \leq &\mathbb{1}_{S_k^i}\mathbb{1}_{\bar J_k^i}\hat c\bm{\Delta}_{k+1}   +\mathbb{1}_{S_k^i}\mathbb{1}_{\bar J_k^i}\frac{\Delta F_{k,0}^i + \Delta F_{k,+}^i+ \Delta F_{k,-}^i}{\bm{\delta}_{k}}.
\end{eqnarray*}

Whereas, when $\mathbb{1}_{\bar S_k^i}=1$, using \eqref{case_b_nuovo} and \eqref{bound_grad_bad_casei}, we get

\begin{eqnarray*}
 \mathbb{1}_{\bar S_k^i}\mathbb{1}_{I_k^i}|\nabla f(\bm{x}_k)^T e_i| & \leq & \mathbb{1}_{\bar S_k^i}\mathbb{1}_{I_k^i}\hat c\bm{\Delta}_{k+1},\\
\mathbb{1}_{\bar S_k^i}\mathbb{1}_{\bar I_k^i}|\nabla f(\bm{x}_k)^T e_i | &\leq& \mathbb{1}_{\bar S_k^i}\mathbb{1}_{\bar I_k^i}\hat c \bm{\Delta}_{k+1} 
+ \mathbb{1}_{\bar S_k^i}\mathbb{1}_{\bar I_k^i}\dfrac{\Delta G_{k,0}^i + \Delta G_{k,+}^i + \Delta G_{k,-}^i}{\bm{\delta}_k}.
\end{eqnarray*}

% Hence, when $x_k\neq x_{k+1}$, using \eqref{LAM_utile1} and \eqref{LAM_utile1_bad} (or \eqref{bound_grad_bad_casei}), we can write
% \begin{equation}\label{bound_grad_success}
% \begin{split}
%  & \mathbb{1}_{J_k}|\nabla f(x_k)^\top e_i| \leq\mathbb{1}_{J_k}\left({c\eps_f(\gamma+2)+L(\sqrt{n}+1)}\right)\Delta_{k+1},\\
%  & \mathbb{1}_{\bar J_k}|\nabla f(x_k)^\top e_i| \leq \mathbb{1}_{\bar J_k}\left(\left({c\eps_f\gamma+L(\sqrt{n}+1)} \right)\Delta_{k+1} +\dfrac{\Delta F_k^i+\bar\Delta F_k^i}{\delta_k}\right),
%  \end{split}
% \end{equation}
% whereas, when $x_k=x_{k+1}$, using \eqref{bound_grad_good_caseii} and \eqref{bound_grad_bad_caseii}, we have
% \begin{equation}\label{bound_grad_failure}
% \begin{split}
% & \mathbb{1}_{J_k}|\nabla f(x_k)^\top e_i| \leq \mathbb{1}_{J_k}\frac{c\eps_f(\gamma+2) + L}{\theta}\Delta_{k+1},\\
% & \mathbb{1}_{\bar J_k}|\nabla f(x_k)^\top e_i| \leq \mathbb{1}_{\bar J_k}\left(\frac{c\eps_f\gamma+L}{\theta} \Delta_{k+1} +\dfrac{\Delta F_k^i+\bar\Delta F_k^i}{\delta_k}\right).
% \end{split}
% \end{equation}
Then, denoting $\nabla_i f(\bm{x}_k) = \nabla f(\bm{x}_k)^\top e_i$, ${\cal F}_k^{-1} = {\cal F}_{k-1,\ell_k-1}$, $\Delta F_k^i = \Delta F_{k,0}^i+\Delta F_{k,+}^i + \Delta F_{k,-}^i$ and $\Delta G_k^i = \Delta G_{k,0}^i+\Delta G_{k,+}^i + \Delta G_{k,-}^i$, we can write 
\[
\begin{split}
    \mathbb{E}[|\nabla_i f(\bm{x}_k)|^2] &= \mathbb{E}\left[\mathbb{1}_{S^i_k}(\mathbb{1}_{\bar J^i_k}+ \mathbb{1}_{J^i_k})|\nabla_i f(\bm{x}_k)|^2 + \mathbb{1}_{\bar S^i_k}(\mathbb{1}_{\bar I^i_k}+ \mathbb{1}_{I^i_k})|\nabla_i f(\bm{x}_k)|^2\right]
    \end{split}
\]

\[
\begin{split}
& \mathbb{E}\left[\mathbb{1}_{S^i_k}(\mathbb{1}_{\bar J^i_k}+ \mathbb{1}_{J^i_k})|\nabla_i f(\bm{x}_k)|^2 \right]= \mathbb{E}\left[\mathbb{E}[\mathbb{1}_{S^i_k}\mathbb{1}_{\bar J^i_k}|\nabla_i f(\bm{x}_k)|^2\Big|{\cal F}_k^{-1}]\right] + \\
&\qquad {\mathbb{E}\left[\mathbb{E}[\mathbb{1}_{S^i_k}\mathbb{1}_{J^i_k}|\nabla_i f(\bm{x}_k)|^2\Big|{\cal F}_k^{-1}]\right]}\\
    &\leq \mathbb{E}\left[\mathbb{E}\left[\left(\hat c\Delta_{k+1} + \dfrac{\Delta F_k^i}{\bm{\delta}_k}\right)^2\Big|{\cal F}_k^{-1}\right]\right] + \mathbb{E}[\mathbb{E}[\hat c^2\Delta_{k+1}^2|{\cal F}_k^{-1}]]\\
    &\leq \mathbb{E}\left[\mathbb{E}\left[2\hat c^2\Delta_{k+1}^2|{\cal F}_k^{-1}\right] + \mathbb{E}\left[ 2\dfrac{(\Delta F_k^i)^2}{\bm{\delta}_k^2}|{\cal F}_k^{-1}\right]\right] + \mathbb{E}[\mathbb{E}[\hat c^2\Delta_{k+1}^2|{\cal F}_k^{-1}]]\\
    &= 3\mathbb{E}\left[\mathbb{E}\left[\hat c^2\Delta_{k+1}^2|{\cal F}_k^{-1}\right]\right] + \mathbb{E}\left[\mathbb{E}\left[ 2\dfrac{(\Delta F_k^i)^2}{\bm{\delta}_k^2}|{\cal F}_k^{-1}\right]\right].%\\
\end{split}
\]
Furthermore,
\[
\begin{split}
& \mathbb{E}\left[\mathbb{1}_{\bar S^i_k}(\mathbb{1}_{\bar I^i_k}+ \mathbb{1}_{I^i_k})|\nabla_i f(\bm{x}_k)|^2 \right]= \mathbb{E}\left[\mathbb{E}[\mathbb{1}_{\bar S^i_k}\mathbb{1}_{\bar I^i_k}|\nabla_i f(\bm{x}_k)|^2\Big|{\cal F}_k^{-1}]\right] + \\
&\qquad {\mathbb{E}\left[\mathbb{E}[\mathbb{1}_{\bar S^i_k}\mathbb{1}_{I^i_k}|\nabla_i f(\bm{x}_k)|^2\Big|{\cal F}_k^{-1}]\right]}\\
    &\leq \mathbb{E}\left[\mathbb{E}\left[\left(\hat c\Delta_{k+1} + \dfrac{\Delta G_k^i}{\bm{\delta}_k}\right)^2\Big|{\cal F}_k^{-1}\right]\right] + \mathbb{E}[\mathbb{E}[\hat c^2\Delta_{k+1}^2|{\cal F}_k^{-1}]]\\
    &\leq \mathbb{E}\left[\mathbb{E}\left[2\hat c^2\Delta_{k+1}^2|{\cal F}_k^{-1}\right] + \mathbb{E}\left[ 2\dfrac{(\Delta G_k^i)^2}{\bm{\delta}_k^2}|{\cal F}_k^{-1}\right]\right] + \mathbb{E}[\mathbb{E}[\hat c^2\Delta_{k+1}^2|{\cal F}_k^{-1}]]\\
    &= 3\mathbb{E}\left[\mathbb{E}\left[\hat c^2\Delta_{k+1}^2|{\cal F}_k^{-1}\right]\right] + \mathbb{E}\left[\mathbb{E}\left[ 2\dfrac{(\Delta G_k^i)^2}{\bm{\delta}_k^2}|{\cal F}_k^{-1}\right]\right].%\\
\end{split}
\]

% Using the Cauchy-Schwarz inequality, we can write
% \[
% \mathbb{E}\left[ \mathbb{1}_{\bar J_k^i}\dfrac{\Delta F_k}{\delta_k}|{\cal F}_k^{-1}\right]\leq \mathbb{E}\left[\mathbb{1}_{\bar J_k^i}|{\cal F}_k^{-1}\right]^{1/2}\mathbb{E}\left[\dfrac{(\Delta F_k)^2}{\delta_k^2}|{\cal F}_k^{-1}\right]^{1/2} 
% \]
%and 
Since $\bm{\delta}_k$ is measurable when conditioning on ${\cal F}_k^{-1}$, we have
% \[
% \begin{split}
% \mathbb{E}\left[ \mathbb{1}_{\bar J_k^i}\dfrac{\Delta F_k}{\delta_k}|{\cal F}_k^{-1}\right]&\leq \mathbb{E}\left[\mathbb{1}_{\bar J_k^i}|{\cal F}_k^{-1}\right]^{1/2}\dfrac{1}{\delta_k}\mathbb{E}\left[{(\Delta F_k)^2}|{\cal F}_k^{-1}\right]^{1/2} 
% \leq \dfrac{\sqrt{1-\beta^3}}{\delta_k}\mathbb{E}\left[{(\Delta F_k)^2}|{\cal F}_k^{-1}\right]^{1/2}\\
% &\leq 
% \dfrac{\sqrt{1-\beta^3}}{\delta_k}\sqrt{3}\eps_f\sqrt{1-\beta}\delta_k^2 = \sqrt{3}\eps_f{\sqrt{1-\beta^3}}\sqrt{1-\beta}\delta_k. 
% \end{split}
% \]
\[
\begin{split}
\mathbb{E}\left[\dfrac{(\Delta F_k^i)^2}{\bm{\delta}_k^2}|{\cal F}_k^{-1}\right]&\leq \dfrac{1}{\bm{\delta}_k^2}\mathbb{E}\left[{(\Delta F_k^i)^2}|{\cal F}_k^{-1}\right] 
%\\&
\leq 
\dfrac{1}{\bm{\delta}_k^2}3c^2\eps_f^2{(1-\beta)}\bm{\delta}_k^4 = 3c^2\eps_f^2{(1-\beta)}\bm{\delta}_k^2\\
&\leq 3c^2\eps_f^2{(1-\beta)}\bm{\Delta}_k^2,\\ 
\mathbb{E}\left[\dfrac{(\Delta G_k^i)^2}{\bm{\delta}_k^2}|{\cal G}_k^{-1}\right]&\leq \dfrac{1}{\bm{\delta}_k^2}\mathbb{E}\left[{(\Delta G_k^i)^2}|{\cal F}_k^{-1}\right] 
%\\&
\leq 
\dfrac{1}{\bm{\delta}_k^2}3c^2\eps_f^2{(1-\beta)}\bm{\delta}_k^4 = 3c^2\eps_f^2{(1-\beta)}\bm{\delta}_k^2\\
&\leq 3c^2\eps_f^2{(1-\beta)}\bm{\Delta}_k^2, 
\end{split}
\]
where the second inequality follows by Assumption \ref{ass:betaprob} and recalling the definition of  $\Delta F_k^i$. %$= |F(\bm{y}_k^i)-f(\bm{y}_k^i)| + |F(\bm{y}_k^i-\bm{\bar\alpha}_k^ie_i)-f(\bm{y}_k^i-\bm{\bar\alpha}_k^ie_i)| + |F(\bm{y}_k^i+\bm{\bar\alpha}_k^ie_i)-f(\bm{y}_k^i+\bm{\bar\alpha}_k^ie_i)|$.

Hence, we obtain
% \[
% \mathbb{E}[|\nabla_i f(x_k)|] \leq (\tilde c+\hat c)\mathbb{E}[\Delta_{k+1}] + 6\eps_f\sqrt{(1-\beta^3)(1-\beta)}\mathbb{E}[\delta_k]
% \]
\[
\mathbb{E}[|\nabla_i f(\bm{x}_k)|^2] \leq 6\hat c^2\mathbb{E}[\bm{\Delta}_{k+1}^2] + 12c^2\eps_f^2{(1-\beta)}\mathbb{E}[\bm{\Delta}_k^2]
\]
which concludes the proof.$\hfill\Box$

\par\medskip

In the following corollary we finally bound the expected value of the square norm  of the gradient  of $f$  by the expected value of $\Delta_{k+1}^2$.

\begin{corollary}\label{bound_grad_coro}
Under the assumptions of Proposition~\ref{bound_grad_prop}, it results
\begin{eqnarray}
%    \mathbb{E}[\|\nabla f(x_k)\|] & \leq & \sqrt{n}\check c\mathbb{E}\left[\Delta_{k+1}\right] + \dfrac{6\eps_f \sqrt{n(1-\beta^3)(1-\beta)}}{\theta}\mathbb{E}[\Delta_{k+1}],\label{bound_norm_grad1}\\
{\mathbb{E}[\|\nabla f(\bm{x}_k)\|^2] }& {\leq }&{6n\hat c^2\mathbb{E}[\bm{\Delta}_{k+1}^2] + \dfrac{12n c^2\eps_f^2(1-\beta)}{\theta^2}\mathbb{E}[\bm{\Delta}_{k+1}^2].\label{bound_norm1_grad1} }
\end{eqnarray}
%{\bf\color{red}dimostrare la relazione in RED}
\end{corollary}
{\bf Proof}. From the instructions of the Algorithm, for all $k$, it results
\[
{\theta\bm{\Delta}_k} \leq \bm{\Delta}_{k+1}.
\]
Then, the proof easily follows from Proposition~\ref{bound_grad_prop}. $\hfill\Box$

% \begin{theorem}
% For the proposed algorithm and under Assumption \ref{ass:betaprob}, we have that
% \[
% \lim_{k\to\infty}\mathbb{E}[\|\nabla f(x_k)\|] = 0.
% \]
% \end{theorem}
% {\bf Proof}.
% The result easily follows by recalling~\eqref{bound_norm_grad} in Proposition \ref{bound_grad_prop} and Theorem \ref{teo:limexpectdelta_tozero}.
% $\hfill\Box$
\par\medskip\noindent
{Finally it is possible to show that the norm of the gradient converges to zero almost surely.}
\begin{theorem}
Let $\beta$ be chosen as in Proposition \ref{phi_decrease}, i.e. 
such that
\[
\frac{\beta^2}{1-\beta^2} > \frac{2\nu}{\min\{{\nu}(\gamma-2) \eta^2,(1-\nu)(1-\theta^2)\}}.
\]
Then, for the proposed algorithm and under Assumption \ref{ass:betaprob}, we have that
\[
\lim_{k\to\infty}\|\nabla f(\bm{x}_k)\| = 0,\quad\text{almost surely}.
\]
\end{theorem}
{\bf Proof}. 
For any $\eps > 0$, we have, for all $k$, 
\[
\mathbb{P}[\|\nabla f(\bm{x}_k)\| > {\eps}] = \mathbb{P}[\|\nabla f(\bm{x}_k)\|^2 > {\eps^2}]. 
\]
By the Markov inequality, we have that
\[
\mathbb{P}[\|\nabla f(\bm{x}_k)\|^2 > {\eps^2}] \leq \frac{\mathbb{E}[\|\nabla f(\bm{x}_k)\|^2]}{\eps^2}
\]
which, from \eqref{bound_norm1_grad1}, gives
\[
\mathbb{P}[\|\nabla f(\bm{x}_k)\|^2 > {\eps^2}] \leq \frac{C}{\eps^2}\mathbb{E}[\bm{\Delta}_{k+1}^2],
\]
where 
\[
C = 6n\hat c_1 + \dfrac{12n c^2\eps_f^2(1-\beta)}{\theta^2}. 
\]
Therefore, summing up for all $k$ we have
\[
\sum_{k=0}^\infty\mathbb{P}[\|\nabla f(\bm{x}_k)\|^2 > {\eps^2}] \leq \frac{C}{\eps^2}\sum_{k=0}^\infty\mathbb{E}[\bm{\Delta}_{k+1}^2].
\]
From Theorem \ref{teo:limexpectdelta_tozero} we know that the summation on the right hand side is finite, hence it is finite also the summation on the left hand side. Hence, we can conclude that 
\[
\mathbb{P}[(\|\nabla f(\bm{x}_k)\|^2 > {\eps^2})\ i.o.] = 0,
\]
thus concluding the proof. $\hfill\Box$

\section{Worst-case complexity result for SDFL}
In this section we derive a worst-case complexity result for algorithm SDFL. In particular, under Assumption \ref{ass:betaprob}, we show that, for any $\eps > 0$, the total number of iterations for which $\mathbb{E}[\|\nabla f(\bm{x}_k)\|]> \eps$ is upper bounded by ${\cal O}(n^2\eps^{-2}/\beta^2)$. Note that, if we denote by $k_\eps$ as the first iteration such that $\mathbb{E}[\|\nabla f(\bm{x}_{k_\eps})\|]\leq \eps$, we also have that $k_\eps\leq {\cal O}(n^2\eps^{-2}/\beta^2)$.  
\begin{proposition}
Let $\beta$ be chosen as in Proposition \ref{phi_decrease}, i.e. 
such that
\[
\frac{\beta^2}{1-\beta^2} > \frac{2\nu}{\min\{{\nu}(\gamma-2) \eta^2,(1-\nu)(1-\theta^2)\}}.
\]
For algorithm SDFL and for any given $\eps>0$, let us consider the set of iterations   
\[
K_\varepsilon = \{k:\ \mathbb{E}[\|\nabla f(\bm{x}_k)\|] > \varepsilon\}=\{k:\ \mathbb{E}[\|\nabla f(\bm{x}_k)\|]^2 > \varepsilon^2\}.
\]
Then, under Assumption \ref{ass:betaprob}, it results
\[
|K_\eps| \leq {\cal O}\left(\frac{n^2\varepsilon^{-2}}{\beta^2}\right).
\]
\end{proposition}
{\bf Proof}. By the Jensen inequality, we have that
\[
\mathbb{E}[\|\nabla f(\bm{x}_k)\|]^2 \leq \mathbb{E}[\|\nabla f(\bm{x}_k)\|^2]
\]
so that we can write, for any $k\in K_\varepsilon$,
\[
 \varepsilon^2 < \mathbb{E}[\|\nabla f(\bm{x}_k)\|^2].
\]

% which, recalling~\eqref{bound_norm1_grad1}, allows us to write%Recalling~\eqref{bound_norm_grad1}, we have that
% \[
% K_\varepsilon \subseteq \bar K_\varepsilon = %\left\{k:\ \check c_3\mathbb{E}\left[\Delta_{k+1}\right]> \varepsilon\right\} = 
% \left\{k:\ \check c_3\mathbb{E}\left[\Delta_{k+1}^2\right]> \varepsilon^2\right\}.
% \]
Now, recalling Proposition~\ref{phi_decrease}, we have, for all $k$,
\[
\mathbb{E}[\Phi_{k+1} - \Phi_k | {\cal F}_{k-1,\ell_k-1}] \leq -\rho\bm{\Delta}_k^2
\]
where 
\[
\rho = \frac{1}{2}\beta^2\min\{{\nu}(\gamma-2) \eta^2,(1-\nu)(1-\theta^2)\},
\]
so that, by taking expectations, we can write
\[
\mathbb{E}[\Phi_{k+1} - \Phi_k] \leq -\rho\mathbb{E}[\bm{\Delta}_k^2].
\]
Now, by summing up over iterations from $0$ to $N$, we have
\[
\rho\sum_{k=0}^N\mathbb{E}[\bm{\Delta}_k^2] \leq \sum_{k=0}^N\mathbb{E}[\Phi_{k}-\Phi_{k+1}] = \mathbb{E}[\Phi_0] - \mathbb{E}[\Phi_{N+1}]\leq \mathbb{E}[\Phi_0]. 
\]
Then, taking the limit for $N\to\infty$,
\[
%\rho|K_\varepsilon|\dfrac{\varepsilon^2}{\check c_3}\leq 
\rho\sum_{k\in K_\varepsilon}\mathbb{E}[\bm{\Delta}_k^2]\leq\rho\sum_{k=0}^\infty\mathbb{E}[\bm{\Delta}_k^2] \leq \mathbb{E}[\Phi_0].
\]
Then, recalling~\eqref{bound_norm1_grad1}, we can write
\[
\rho|K_\varepsilon|\dfrac{\varepsilon^2}{\check c_3}\leq
\rho\sum_{k\in K_\varepsilon}\mathbb{E}[\bm{\Delta}_k^2]
\leq\mathbb{E}[\Phi_0].
\]

Hence, we obtain
\[
\begin{split}
|K_\varepsilon|&\leq \dfrac{\check c_3\mathbb{E}[\Phi_0]}{\rho\varepsilon^2} = \dfrac{\check c_3 \Phi_0}{\rho\varepsilon^2} = 2\dfrac{3n\hat c^2 + \dfrac{6n c^2\eps_f^2(1-\beta)}{\theta^2}}{\beta^2\min\{{\nu}(\gamma-2) \eta^2,(1-\nu)(1-\theta^2)\}}\Phi_0\\
&\leq 2\dfrac{6n\hat c^2 + \dfrac{12n c^2\eps_f^2}{\theta^2}}{\beta^2\min\{{\nu}(\gamma-2) \eta^2,(1-\nu)(1-\theta^2)\}}\Phi_0= {\cal O}\left(\frac{n^2\varepsilon^{-2}}{\beta^2}\right),
\end{split}
\]
thus concluding the proof. $\hfill\Box$

\section{Conclusions}
In this paper we have studied the convergence properties and the worst case iteration complexity of a derivative-free algorithm based on extrapolation techniques when applied to a stochastic problem, i.e. problem \eqref{theoretical_prob}. The analysis of our algorithm SDFL is inspired by the recent papers \cite{dzahini2022expected,audet2021stochastic}. However, the probabilistic convergence properties we derived for our method are somewhat stronger than those proved in \cite{dzahini2022expected,audet2021stochastic}. More in particular, 
% \begin{enumerate}
%     \item in \cite{dzahini2022expected}, it is proved that
%     \[
%     \liminf_{k\to\infty}\|\nabla f(\bm{x}_k)\| = 0\qquad\text{almost surely}
%     \]
%     \item in \cite{audet2021stochastic}, it is proved that the Clarke generalized derivative of $f$ at a refined point along any corresponding refined direction is almost surely non-negative. 
% \end{enumerate}
For our SDFL algorithm We managed to prove that
\[
\lim_{k\to\infty}\|\nabla f(\bm{x_k})\| = 0\qquad\text{almost surely},
\]
i.e. that every limit point of the sequence of stochastic iterates is almost surely stationary.

Concerning the worst case complexity, our result is similar but different than the complexity proved in \cite{dzahini2022expected}. Specifically, in \cite{dzahini2022expected} it is proved that $\mathbb{E}[T_\eps^*]\leq{\cal O}(\eps^{-p/\min\{p-1,1\}}/(2\beta-1))$ where $T^*_\eps$ is the so-called {\it stopping time}, i.e.
\[
T^*_\eps = \inf\{k:\ \|\nabla f(\bm{x}_k)\|\leq\eps\}
\]
for a given $\epsilon>0$. Instead, for algorithm SDFL, we managed to prove that $|K_\eps| \leq {\cal O}(\eps^{-2}/\beta^2)$ where $|K_\eps|$ is the total number of iterations where the expected value of the gradient norm is above $\eps$, i.e. $\mathbb{E}[\|\nabla f(\bm{x}_k)\|]>\eps$. It is worth noticing that the role played by $\beta$ in \cite{dzahini2022expected}, that is the probability of the intersection event that two function values are accurate, is the same played by our $\beta^2$ coefficient. Indeed, we define $\beta$ as the probability that a single function value is accurate so that $\beta^2$ is the probability that two function values are contemporaneously  accurate.

Furthermore, we note that our complexity result also gives a bound on the first iteration such that the expected value of the gradient norm is below a prefixed tolerance. However, such a bound, though different from those in \cite{dzahini2022expected}, is worst in terms of the constants multiplying the $\eps^{-2}$. Indeed, we have a $n^2$ instead of an $n$. Considering the recent paper on the complexity of linesearch derivative-free methods, it is possible to improve the bound to ${\cal O}(n\eps^{-2})$. However, this further analysis is not trivial in the stochastic context and would surely be the subject of future work.

%\bibliography{noisybiblio}

\end{document}